\documentclass{article}

\usepackage[margin=0.94in]{geometry}
\usepackage{graphicx,parskip} % Required for inserting images
\usepackage{amsmath,amsfonts,amssymb,mathdots}
\usepackage{amsthm}
\usepackage{xcolor}

\usepackage{ytableau}
\usepackage{url}
\usepackage{hyperref}

\newtheorem{theorem}{Theorem}[section]
\newtheorem{prop}[theorem]{Proposition}

\newtheorem{lemma}[theorem]{Lemma}
\newtheorem{corollary}[theorem]{Corollary}
\newtheorem{remark}[theorem]{{\sc Remark}}
\newtheorem{example}[theorem]{Example}

\theoremstyle{definition}
\newtheorem{definition}[theorem]{Definition}

\newenvironment{keywords}{ {\bf Keywords:}}

\newenvironment{AMS}{ {\bf AMS:}}

\newcommand{\rank}{\mathrm{rank}}
\newcommand{\Res}{\mathrm{Res}}
\newcommand{\ideal}[1]{\langle #1 \rangle}
\newcommand{\st}{~\mathrm{s.t.}~}
\newcommand{\LM}{\mathrm{LM}}

\newcommand{\FF}{\mathbb{F}}

\newcommand{\NN}{\mathbb{N}}

\newcommand{\BB}{\mathcal{B}}

\definecolor{vannihomepage}{rgb}{0.91, 0.69, 1.0}
\definecolor{homepagedark}{rgb}{0.85, 0.48, 1.0}
\definecolor{cutepink}{rgb}{0.76, 0.33, 0.46}
\definecolor{skyblue}{rgb}{0.22, 0.57, 0.89}
\definecolor{myviolet}{rgb}{0.21, 0.0, 0.85}

\title{The Smith form of Sylvester and B\'{e}zout matrices\\ for zero-dimensional ideals}
\author{Etna Lindy\footnotemark[1] \and Vanni Noferini\footnotemark[2]}

\begin{document}

\maketitle

\renewcommand{\thefootnote}
{\fnsymbol{footnote}}
\footnotetext[1]{Aalto University, Department of Mathematics and Systems Analysis, P.O. Box 11100, FI-00076, Aalto, Finland. Supported by the Finnish Ministry of Education and Culture’s Pilot for Doctoral Programmes (Pilot project Mathematics of Sensing, Imaging and Modelling). Email: etna.lindy@aalto.fi}
\footnotetext[2]{Corresponding author. Aalto University, Department of Mathematics and Systems Analysis, P.O. Box 11100, FI-00076, Aalto, Finland. Supported by a Research Council of Finland grant (decision number 370932). Email: vanni.noferini@aalto.fi}
\renewcommand{\thefootnote}{\arabic{footnote}}

\begin{abstract}
 Let $\mathbb{K}$ be a field and let $f,g \in \mathbb{K}[x,y]$ be such that the ideal $\ideal{f,g}$ is zero-dimensional. We study the Sylvester and B\'{e}zout resultant polynomial matrices, built by interpreting $f$ and $g$ as univariate polynomials in $x$ with coefficients in $\mathbb{K}[y]$. We characterize their Smith forms over $\mathbb{K}[y]$ in terms of the dual spaces of differential operators, that were defined and studied by H. M. M\"{o}ller et al. In particular, if $\mathbb{K}$ is algebraically closed we show that, if the leading coefficients of $f$ and $g$ are coprime over $\mathbb{K}[y]$, then the partial multiplicities of the Sylvester and B\'{e}zout resultant matrices coincide with certain integers, that we call M\"{o}ller indices. These indices are uniquely determined by $\ideal{f,g}$, and can be easily computed from a Gauss basis, as defined in [M. G. Marinari, H. M. M\"{o}ller, T. Mora, Trans. Amer. Math. Soc. 348(8):3283--3321, 1996], of the dual spaces. We then generalize this result to the case of common factors in the leading coefficients, which correspond to intersections at $x=\infty$, again describing all the invariant factors of Sylvester and B\'{e}zout resultant matrices. As a corollary, this fully characterizes the algebraic multiplicity of all the roots of the resultant $\Res_x(f,g) \in \mathbb{K}[y]$ in terms of the intersection multiplicities for $f$ and $g$, including those arising from infinite intersections. We discuss both algebraic and computational implications of our results.
\end{abstract}
\begin{keywords}
Resultant matrix, invariant factor, Smith form, zero-dimensional ideal, dual space, Gauss basis, M\"{o}ller index
\end{keywords}
\begin{AMS}
13P15, 68W30, 14Q99, 15A18, 15A22
\end{AMS}

\section{Introduction}

In elimination theory, a central role is played by classical resultant matrices such as those of Sylvester and B\'{e}zout. Suppose for a start that $f$ and $g$ are nonconstant univariate polynomials in $x$ with coefficients in an integral domain $R$; then, the associated Sylvester and B\'{e}zout matrices are square matrices with entries in $R$, and both of them are singular if and only if $f$ and $g$ have a common root in the algebraic closure of $R$. In fact, the \emph{resultant} of $f$ and $g$ is often defined to be the determinant of the Sylvester matrix \cite{iva,usingag}; instead, the determinant of the B\'{e}zout matrix is the resultant times a certain ``extraneous factor", equal to the leading coefficient of one of the two polynomials raised to the difference of the degrees of the two polynomials.

If we take $R=\FF[y]$, where $\FF$ is an algebraically closed\footnote{The assumption that $\FF$ is algebraically closed significantly simplifies the exposition, but it can be weakened. Namely, if $f,g \in \mathbb{K}[x,y]$, we can embed the arbitrary field $\mathbb{K}$ in its algebraic closure $\FF=\overline{\mathbb{K}}$, and note that this field extension leaves the monic invariant factors of the associated resultant matrices unchanged.} field, we can leverage the notion of resultant as a tool to solve the bivariate polynomial system $f(x,y)=g(x,y)=0$. When the ideal $\ideal{f,g}$ is zero-dimensional, the Sylvester and B\'{e}zout resultant matrices are commonly used as a tool to solve polynomial equations \cite{BM,GT,MD,NNT,NT,petit,sommese,villard}; for this reason, we focus on zero-dimensional ideals in this paper. Indeed, it is clear that if $(x_0,y_0) \in \FF^2$ is a solution to the system, then the resultant of $f$ and $g$, which is a polynomial in $\FF[y]$, must vanish at $y=y_0$; more details can be found, for example, in \cite[Chapter 3]{iva} or \cite[Chapters 3 and 7]{usingag}. The reverse implication is subtler. If the leading coefficients of $f$ and $g$ do not share any non-unital common factor in $\FF[y]$, then, for every root $y_0 \in \FF$ of the resultant, there indeed exists $x_0 \in \FF$ such that $f(x_0,y_0)=g(x_0,y_0)=0$ \cite[Section 3.6, Corollary 7]{iva}. Without assuming coprime leading coefficients, however, it may happen that $y_0$ is a root of the resultant, but still $\begin{bmatrix}
    f(x_0,y_0)\\
    g(x_0,y_0)
\end{bmatrix} \neq 0 \ \forall \ x_0 \in \FF$. One may express this fact by saying that \cite{NNT} the roots of the resultant correspond to either the $y$-components of the finite solutions of the polynomial system of equations, or the $y$-components of ``roots at infinity" (in the $x$-component). In the computational algebraic geometry literature, this fact may be considered folklore in the sense that it is often acknowledged in passing to constitute a complication \cite{iva,usingag,NNT,villard}, but a rigorous structural characterization of the relation between these ``roots at infinity" and the roots of the resultant (or, even more deeply, the roots of each invariant factor of the resultant matrices) seems to be lacking; one major contribution of the present paper is to explicitly provide a precise analysis, by deriving the full Smith form of the Sylvester and B\'{e}zout resultant matrices in the case where $\langle f,g \rangle$ is zero dimensional. In particular, such Smith forms yield all the partial multiplicities of the elementary divisor $y-y_0 \in \FF[y]$, where $y_0 \in \FF$ is either the $y$-component of a finite intersection $(x_0,y_0) \in \FF^2$ or of an intersection at infinity $(\infty,y_0)$ (Definition \ref{def:projectiveintersection}); this has computational consequences, because when $\mathbb{F} \subseteq \mathbb{C}$ the partial multiplicities determine whether $y_0$ is ill conditioned or not as an eigenvalue of the resultant matrix.

Related results that can be found in the literature include:
\begin{itemize}
    \item \cite[Theorem 4]{lazard}, where D. Lazard described (implicitly, via a special Hermite form) the Smith form of the Sylvester matrix in terms of Gr\"{o}bner bases. However, the focus of \cite{lazard} is elsewhere, and a structural analysis of the invariant factors is not explicitly emphasized. Moreover, \cite[Theorem 4]{lazard} assumed that the leading coefficients of $f$ and $g$ are coprime, and \cite[Remark 7]{lazard} explicitly posed as an open problem the analysis of the general case.
    \item \cite[Theorem 1.1]{bda}, where L. Bus\'{e} and C. D'Andrea computed the invariant factors of certain (specific) Sylvester and B\'{e}zout bivariate homogenous resultant matrices arising from birational parameterizations of rational plane curves, describing them in terms of the singular point of such curves. The characterization of \cite{bda} is explicit, but restricted to a special class of polynomials.
    \item Several sources, such as for instance \cite[Lemma 2.2]{villard}, identified the last invariant factor of the Sylvester matrix as a generator of the elimination ideal $\langle f,g \rangle \cap \FF[y]$.
\end{itemize}

In this paper, we generalize and complete this picture. First, under the assumption of coprime leading coefficients of $f$ and $g$, we provide a dual characterization of the invariant factors in terms of the \emph{dual spaces}  of linear functionals \cite{alonso,3M91,3M,3M96,moller,ms} vanishing on the ideal $\langle f,g \rangle$; see also \cite{moupan} for a thorough review on dual spaces and structured matrices. We connect the local dual spaces associated with  the ideal $\langle f,g \rangle$ to the the Smith form of the resultant matrices by means of the theory of root vectors \cite{rp1,rp3,GLR82,rp2}. If we fix one point $(x_0,y_0)$ on the variety of $\ideal{f,g}$ and we endow $\FF[x,y]$ with the lexicographic monomial order such that $x<y$, then the shape of the local dual space at $(x_0,y_0)$ is described by certain combinatorial invariants, that we propose to call the \emph{M\"{o}ller indices with respect to $y$} of $\langle f,g \rangle$ at $(x_0,y_0)$ (Definition \ref{def:mollerindices}). These indices are uniquely determined by $\langle f, g \rangle$ (Remark \ref{rem:mollerindices}) and can be easily computed from any Gauss basis \cite{3M,3M96} (Theorem \ref{thm:mollerindequiv}). We prove in Theorem \ref{thm:mainresultnoinf} that, for all $y_0 \in \FF$, the partial multiplicities of  the elementary divisor $y-y_0$ in the Sylvester matrix of $f$ and $g$ coincide with the M\"{o}ller indices with respect to $y$ of $\langle f,g \rangle$ at all points of the form $(x_i,y_0)$, $x_i \in \FF$, in the variety of the ideal $\ideal{f,g}$. This is a dual description with respect to the approach of \cite[Theorem 4]{lazard}, based on Gr\"{o}bner bases. Furthermore, unlike \cite{lazard} we also give an explicit construction of root vectors based on a Gauss basis of the local dual space; root vectors go beyond the partial multiplicities, providing information on the unimodular matrices in a Smith form decomposition \cite{rp3,rp2} of the Sylvester or the B\'{e}zout matrices.

Next, we extend our analysis to also cover infinite roots, i.e., we allow the leading coefficients of $f$ and $g$ to have a nontrivial GCD over $\FF[y]$, while keeping the assumption that $\langle f,g \rangle$ is zero-dimensional. In this situation, and after appropriately defining roots at infinity, the infinite roots also contribute to the Smith form of the resultant matrices, in a completely analogous manner as the finite roots. Thus, in Theorem \ref{thm:mainresult} and Theorem \ref{thm:mainresultbezout}, respectively for Sylvester and B\'{e}zout matrices, we provide a unified treatment of finite and infinite intersections within the same algebraic framework.

Our results have both algebraic and computational implications. In the case where $\langle f,g \rangle$ is zero-dimensional and the leading coefficients (in $x$) of $f$ and $g$ are $\FF[y]$-coprime, then the Smith form of the Sylvester matrix of $f$ and $g$ yields the structure of $\FF[x,y]/\langle f,g \rangle$ as an $\FF[y]$-module (Theorem \ref{thm:wehatealggeom}). We thus describe such module in terms of the M\"{o}ller indices. On the other hand, for $\mathbb{F} \subseteq \mathbb{C}$ and from a numerical linear algebra viewpoint, the partial multiplicities of the resultant matrices determine the H\"{o}lder condition numbers \cite{KPM} of their eigenvalues. Thus, the partial multiplicities, and therefore the M\"{o}ller indices, predict the expected accuracy of the numerical computation of eigenvalues of resultant matrices, which lies at the core of many resultant-based rootfinding methods \cite{GT,NNT}. Namely, even if the intersection multiplicity is $>1$, the eigenvalue corresponding to the $y$-component of the intersection may still be computable numerically with high accuracy if, say, the M\"{o}ller indices are all $1$; and conversely, high M\"{o}ller indices imply that one should have low expectations for the reliability of floating point computations of the associated eigenvalues. In some applications, including cryptography \cite{petit} and number theory \cite{cohen}, polynomial rootfinding also arises over other fields than the complex numbers, and our results may possibly benefit algorithmic approaches in these areas as well.

The manuscript is structured as follows. Section \ref{sec:bg} recalls the necessary preliminary material, including dual spaces, resultant matrices, Smith forms, and root vectors. In Section \ref{sec:indices} we define M\"{o}ller indices and show how they can be computed from a Gauss basis \cite{3M,3M96}. The next sections are devoted to stating and proving our main results for the Sylvester (Section \ref{sec:sylvestersmith}) and B\'{e}zout matrix (Section \ref{sec:bezoutsmith}). Finally, in Section \ref{sec:numerical} we make some numerical considerations for $\FF \subseteq \mathbb{C}$.

\section{Background}\label{sec:bg}

Let us begin by setting some basic notation. Throughout $\FF$ will denote an algebraically closed field, and $\FF[x]$ a univariate polynomial ring in the variable $x$ and with coefficients in $\FF$. Similarly, $\FF[x,y]$ denotes bivariate polynomials in $x$ and $y$.
Sometimes we perceive elements of $\FF[x,y]$ as elements of $\FF[y][x]$, that is, univariate polynomials in $x$ with coefficients in $\FF[y]$. %When we wish to enforce that the coefficients to lay in a field, we utilize the extension of $\FF[y]$ to the field of fractions $\FF(y)$, or even the algebraic closure $\PP(y)$. 

Given an integral domain $R$ and $n \in \mathbb{N}$, we denote by $R[x]_n$ the $R$-module of polynomials in $x$ with coefficients in $R$ and degree at most $n$. %Similarly, the denotation $R[x,y]_{m,n}$ denotes the $R$-module of bivariate polynomials in $x,y$ and coefficients in $R$ having separate degree at most $m$ in $x$ and at most $n$ in $y$. 
For every positive integer $k>0$, we will make frequent use of the Vandermonde vector
$ \Lambda_k(x) : = \begin{bmatrix}
    x^{k-1}\\
    \vdots\\
    x^2\\
    x\\
    1
\end{bmatrix} \in R[x]^k.$ We also often employ the $i$-th Hasse derivative (see Subsection \ref{sec:functionals}) of $\Lambda_k(x)$, denoted by

\[ \Lambda_k^{(i)}(x)  = \begin{bmatrix}
   \begin{pmatrix}
       k-1\\
       i
   \end{pmatrix} x^{k-i-1}\\
    \vdots\\
   \begin{pmatrix}
       i+2\\
       i
   \end{pmatrix} x^2\\
   (i+1) x\\
    1\\
    0\\
    \vdots\\
    0
\end{bmatrix} \in R[x]^k.\] Finally, for $x_0 \in R$ the notation $\Lambda_k^{(i)}(x_0)$ simply means the vector of $R^k$ obtained by evaluating $\Lambda_k^{(i)}(x)$ at $x=x_0$, and similarly for other polynomial vectors or matrices.   We also denote the ``coefficient extraction map" by
\[  \Gamma_n : R[x]_{n-1} \rightarrow R^n, \qquad f(x)=\sum_{i=0}^{n-1} f_i x^i \mapsto \Gamma_n(f)= \begin{bmatrix}
    f_{n-1}\\
    \vdots\\
    f_1\\
    f_0
\end{bmatrix}.         \]
Note that the identity $f(x)=\Gamma_n(f)^\top \Lambda_n(x)$ holds for all $f(x) \in R[x]_{n-1}$. 
Just as we match vectors and univariate polynomials, we can identify a bivariate polynomial and a matrix containing its coefficients in the monomial basis, $f(x,y)=\Lambda_m(y)^\top M \Lambda_n(x)$  \cite{NNT2,rp2}. 

    \subsection{Hasse derivatives, dual spaces, and multiplicities}\label{sec:functionals}

In this subsection, we revisit the theory developed in \cite{3M91,3M,3M96,moller,ms}. Here, we specialize it to the case of two variables. We begin by defining the evaluation functional.
\begin{definition}
   Given a point $(x_0,y_0) \in \FF^2$, the {\it evaluation map} is the functional $E_{(x_0,y_0)} : \FF[x,y] \rightarrow \FF$ that satifies  
    $E_{(x_0,y_0)}f := f(x_0,y_0)$ for all $f \in \FF[x,y]$.
\end{definition}

When the value of one coordinate is irrelevant, we may write, e.g., $E_{(\cdot,y_0)}$. For instance, if $g \in \FF[y]$ then $E_{(\cdot,y_0)}g=g(y_0) \in \FF$. Similarly, if $h \in \FF[x,y]$ then $E_{(x_0,\cdot)}h=h(x_0,y) \in \FF[y]$. For every $(i,j) \in \NN^2$, the second step is to define the Hasse derivatives as the $\FF$-linear operators $D_{ij}:\FF[x,y] \rightarrow \FF[x,y]$ satisfying
\[  D_{ij} x^a y^b =  \begin{cases}
 \begin{pmatrix}
     a\\
     i
 \end{pmatrix} \begin{pmatrix}
     b\\
     j
 \end{pmatrix} x^{a-i} y^{b-j} & \ \mathrm{if} \ a \geq i \ \mathrm{and} \ b \geq j;\\
  0 & \ \mathrm{otherwise.}
\end{cases}     \]
We note the useful properties $D_{ij}=D_{0j}D_{i0}=D_{i0}D_{0j}$ and, for all $f,g \in \FF[x,y]$, \[D_{ij}fg=\sum_{k=0}^i \sum_{\ell=0}^j D_{i-k,j-\ell} f D_{k\ell} g.\] When $\FF$ has zero characteristic, we also observe that $D_{ij} = \frac{1}{i!j!} \frac{\partial^i }{\partial x^i} \frac{\partial^j}{\partial y^j}.$
The following property, whose elementary proof is omitted, will be useful\footnote{If $\FF$ has prime characteristic $p$, it is not possible to replace Hasse derivatives with usual formal derivatives in Proposition \ref{prop:primechar}, e.g., every formal derivative of  $y^p$ is $0$ at $y=0$.}. 
\begin{prop}\label{prop:primechar}
    A polynomial $f(y) \in \FF[y]$ has a root of multiplicity $\mu$ at $y_0 \in \FF$, i.e., $f(y)=(y-y_0)^\mu g(y)$ and $g(y_0) \neq 0$, if and only if $E_{(\cdot,y_0)}D_{0 \mu}f(y) \neq 0$ and $E_{(\cdot,y_0)}D_{0i} f(y) = 0$ for all $0 \leq i < \mu$.
\end{prop}
 Let us now define the $\FF$-vector space of Hasse derivatives
$ \mathcal{D} :=  \mathrm{span}_\FF\{ D_{ij} \mid i,j \in \NN \},$ of which $\{D_{ij}\}_{i,j \in \mathbb{N}}$ form a basis.
Two further ingredients are the antiderivatives with respect to either variable, defined as the $\FF$-linear maps $A_x : \mathcal{D} \to \mathcal{D}$ and $A_y: \mathcal{D} \to \mathcal{D}$ such that
\[ A_x(D_{ij}) = \begin{cases}
    D_{i-1,j}, &\text{if } i > 0 \\
    0  &\text{otherwise},
\end{cases} \qquad  A_y(D_{ij}) = \begin{cases}
    D_{i,j-1}, &\text{if } j > 0 \\
    0  &\text{otherwise}.
\end{cases}\]

\begin{remark}
    All the functionals discussed above (derivatives, antiderivatives, and evaluation maps) may also be applied, elementwise, to bivariate polynomial vectors or matrices.
\end{remark}

\begin{definition}\label{def:closedsubspace}
    A subspace $V \subset \mathcal{D}$ is {\it closed} if it is finite dimensional and
    \[  \phi \in V \implies A_x \phi \in V, ~A_y \phi \in V. \]
\end{definition}
It turns out \cite{3M91,3M,3M96,moller,ms} that there is a correspondence between ideals of $\FF[x,y]$ and subspaces of $\mathcal{D}$.
\begin{definition}
Given $f,g \in \FF[x,y]$, let $I = \ideal{f,g}$ and $(x_0,y_0) \in \mathcal{V}(I)$. The dual space $V_{(x_0,y_0)} \subset \mathcal{D}$ of $I$ at $(x_0,y_0)$ is the subspace
    \begin{align*}
        V_{(x_0,y_0)} = \{ \phi \in \mathcal{D} \ \mathrm{s.t.} \ E_{(x_0,y_0)} (\phi p) = 0, ~\forall p \in I  \}.
    \end{align*}
\end{definition}

M. G. Marinari, H. M. M\"{o}ller and T. Mora proved that the dual space $V_{(x_0,y_0)}$ satisfies $\phi \in V \Rightarrow A_x \phi, A_y \phi \in V$ \cite[Proposition 2.4]{3M} and that it is finite dimensional when $I$ is zero-dimensional \cite[Theorem 3.2]{3M96}; hence, when $I$ is zero-dimensional, the dual space at $(x_0,y_0) \in \mathcal{V}(I)$ is a closed subspace. In that case, following \cite{3M,3M96,ms}, we define the {\it multiplicity} of $(x_0,y_0)$ in $I$ as $\mu(x_0,y_0) := \dim_\FF V_{(x_0,y_0)}$. Lemma \ref{lem:multiplicity} below shows that this definition coincides with the classical definitions of intersection multiplicity.

\begin{lemma} \label{lem:multiplicity}
    Let $I = \ideal{f,g}$ be a zero-dimensional ideal of $\FF[x,y]$ as above and let $P=(x_0,y_0) \in \mathcal{V}(I)$. Consider the local ring 
    \[ \mathcal{O}_{P} := \{ p/q \ \mathrm{s.t.} \ p,q \in \FF[x,y], ~q(x_0,y_0) \neq 0 \}. \]
    The multiplicity of the point $P$ in $I$, as defined above, satisfies
    \[ \mu(x_0,y_0) = \dim_\FF \mathcal{O}_{P} / \ideal{f,g}_\mathcal{O}, \]
    where $\ideal{f,g}_\mathcal{O}$ is considered as an ideal of $\mathcal{O}_{P}$.
\end{lemma}
\begin{proof}
    The statement follows by \cite[Exercise 5.4]{hartshorne} and \cite[Theorem 2.6]{3M}.
\end{proof}

\subsection{Basic properties of Sylvester and B\'{e}zout matrices}

Given two polynomials $f,g \in R[x]$, where $R$ is an integral domain, we now recall the definitions of the associated Sylvester and B\'{e}zout matrices, as well as some of their basic properties. Most, if not all, of this subsection's content is classical \cite{fuhrmann}, but we give our own treatment for self-containedness.

\begin{definition}[Sylvester matrix, \cite{iva}]\label{def:smatrix}
    Let $m \geq \deg f$ and $n \geq \deg g$, $m+n \geq 1$, and write $f(x)=\sum_{i=0}^m f_i x^i, g(x)=\sum_{i=0}^n g_i x^i$ (note that we possibly allow $f_m=0$ and/or $g_n=0$). The $(m,n)$-Sylvester matrix of the univariate polynomials $f$ and $g$ is the matrix 
    \[ S_{f,g;m,n}  := \begin{bmatrix}
        f_{m} & \dots & f_0 & 0 & \dots & 0 \\
        0 & f_{m} & \dots & f_0 & & \vdots &\\
        \vdots &  &  \ddots &  & \ddots & 0\\
        0 & \dots & & f_{m} & \dots & f_0 \\
        g_{n} & \dots & g_0 & 0 & \dots & 0 \\
        0 & g_{n} & \dots & g_0 & & \vdots &\\
        \vdots &  &  \ddots &  & \ddots & 0\\
        0 & \dots & & g_{n} & \dots & g_0 \\
    \end{bmatrix} \in R^{(m+n) \times (m+n)}.\]
     For the choice $m=\deg f,n=\deg g$, we will simply refer to the matrix above as the \emph{Sylvester matrix} of $f$ and $g$.
\end{definition}

Unless a possible ambiguity arises, we will avoid the baroque notation $S_{f,g;m,n}$ and simply write $S$. We omit the straightforward proof of Proposition \ref{prop:sproperty} below.

\begin{prop}\label{prop:sproperty}
    Let $S$ be the $(m,n)$-Sylvester matrix of $f,g$, where $m \geq \deg f, n \geq \deg g$.
    \begin{enumerate}
        \item $S$ is the unique element of $R^{(m+n) \times (m+n)}$ satisfying $S \Lambda_{m+n}(t) = \begin{bmatrix}
            f(t) \Lambda_n(t)\\
            g(t) \Lambda_m(t)
        \end{bmatrix} $.
        \item For arbitrary $p(t) \in R[t]_{n-1},q(t)\in R[t]_{m-1}$, write $p(t)=\Gamma_{n}(p)^\top \Lambda_n(t),q(t)=\Gamma_{m}(q)^\top \Lambda_m(t)$ and define  $r(t):=f(t)p(t)+g(t)q(t)=\Gamma_{m+n}(r)^\top \Lambda_{m+n}(t)$. Then, $S^\top \begin{bmatrix}
            \Gamma_{n}(p)\\
            \Gamma_{m}(q)
        \end{bmatrix} = \Gamma_{m+n}(r).$ In other words, $S^\top$ represents (in the monomial basis) the Sylvester operator
        \[ \mathcal{S} : R[t]_{n-1} \oplus R[t]_{m-1} \rightarrow R[t]_{m+n-1}, \qquad (p(t),q(t)) \mapsto f(t) p(t) + g(t) q(t).  \]
    \end{enumerate}
\end{prop}

\begin{definition}[B\'{e}zoutian and B\'{e}zout matrix, \cite{fuhrmann}]\label{def:bezout}
Given the polynomials $f,g \in R[t]$, their B\'{e}zoutian is the bivariate polynomial
\[ \mathcal{B}_{f,g}(x,z) := \frac{f(x)g(z) - f(z)g(x)}{x-z}  \in R[x,z].   \]
   Let $k=\max \{\deg f, \deg g \} \geq 1$. The B\'{e}zout matrix of $f$ and $g$, $B_{f,g} \in \FF^{k \times k}$, is the unique matrix satsfying
    \[ \Lambda_{k}(z)^\top B_{f,g} \Lambda_k(x) = \mathcal{B}_{f,g}(x,z) = \frac{f(x)g(z) - f(z)g(x)}{x-z}. \]
\end{definition}

Again, we will simply write $\mathcal{B}(x,z)$ and $B$ if there is no ambiguity about the polynomial arguments. Although its construction is more involved than the Sylvester matrix, there are various algorithms to compute the B\'{e}zout matrix, of which (to our knowledge) the most efficient is the one based on ``shifted sums" \cite{NNT,NNT2}. Another method, less effective computationally but more useful for theoretical considerations, is via the connection between the Sylvester matrix and the Bezout matrix that we give in Proposition \ref{prop:SFS} below.

\begin{prop}\label{prop:SFS}
    Let $k=\max\{\deg f, \deg g\} \geq 1$ and let $S$ be the $(k,k)$-Sylvester matrix of $f$ and $g$; let $B$ be the B\'{e}zout matrix of the same pair. Moreover, let $F \in R^{k \times k}$ be the flip matrix, $F_{ij}=\begin{cases}
        1 \ \mathrm{if} \ i+j=k+1;\\
        0 \ \mathrm{otherwise}.
    \end{cases}$ Then,
    \[  S^\top \begin{bmatrix}
        0 & F\\
        -F & 0
    \end{bmatrix} S = \begin{bmatrix}
        0 & B\\
        -B & 0
    \end{bmatrix}.     \]
\end{prop}
\begin{proof}
   Following  \cite{rp3,NNT2}, we note that $ M \in R^{2k \times 2k} \mapsto \Lambda_{2k}(z)^\top M \Lambda_{2k}(x) \in R[x,z]_{2k-1,2k-1} $
    is a bijection. Hence, we may equivalently interpret the statement as an identity for bivariate polynomials.    Using Proposition \ref{prop:sproperty},
    \begin{align*}
        S\Lambda_{2k}(x) = \begin{bmatrix}
            f(x)\Lambda_k(x) \\g(x) \Lambda_k(x)
        \end{bmatrix} \quad \mathrm{and} \quad 
        \Lambda_k(z)^\top F \Lambda_k(x) = \sum_{i=0}^{k-1} x^i z^{k-1-i} = \frac{x^k - z^k}{x - z}.        
    \end{align*}
    Therefore,
    \begin{align*}
        &\Lambda_{2k}(z)^\top S^\top \begin{bmatrix}
        0 & F \\ -F & 0
    \end{bmatrix} S \Lambda_{2k}(x) = [ f(z) \Lambda_k(z)^\top , g(z)\Lambda_k(z)^\top ] \begin{bmatrix}
        g(x) F \Lambda_k(x) \\ -f(x) F \Lambda_k(z)
    \end{bmatrix} \\&=  f(z)g(x) \frac{x^k-z^k}{x-z}  - f(x)g(z) \frac{x^k-z^k}{x-z} = \mathcal{B}(x,z) (z^k - x^k).
    \end{align*}
    On the other hand, by Definition \ref{def:bezout},
    \begin{align*}
        &\Lambda_{2k}(z)^\top \begin{bmatrix}
            0 & B \\ -B & 0
        \end{bmatrix} \Lambda_{2k}(x) = \Lambda_{2k}(z)^\top \begin{bmatrix}
            B \Lambda_k(x) \\
            -x^k B\Lambda_k(x)
        \end{bmatrix}\\
        &= z^k \mathcal{B}(x,z) - x^k \mathcal{B}(x,z) = \mathcal{B}(x,z) (z^k - x^k).
    \end{align*}
    This completes the proof.
\end{proof}

The resultant of $f$ and $g$ is defined \cite{iva,usingag} as the determinant of the Sylvester matrix.

\begin{definition}[Resultant, \cite{iva,usingag}]\label{def:res}
    Let $\FF$ be the algebraic closure of the integral domain $R$. For all polynomials $f = \sum_{i=0}^m f_i x^i \in R[x]$ and $g = \sum_{i=0}^n g_i x^i \in R[x]$, having degree at most $m$ and $n$ respectively, there exists a polynomial $\Res \in R[f_0,\dots,f_{m},g_0,\dots,g_{n}]$,  called the \emph{resultant}, such that
    \[ \Res(f_0,\dots,f_{m},g_0,\dots,g_{n}) = 0 \iff \exists~x_0 \in \FF \st f(x_0) = g(x_0) = 0. \]
    Once we fix $f$ and $g$, i.e., once we specify their coefficients, the resultant $\Res(f,g) \in R$ is uniquely determined by $f$ and $g$ up to a constant; we normalize it by setting $\Res(f,g): = \det(S)$ where $S$ is the $(\deg(f),\deg(g))$-Sylvester matrix of $f$ and $g$.
\end{definition}

The determinant of the B\'{e}zout matrix is a close relative of the resultant.
\begin{corollary}\label{cor211}
Let $f,g$ have degrees $m,n$ respectively, and let $k=\max\{m,n\}$.
    If $B \in R^{k \times k}$ is the B\'{e}zout matrix of $f$ and $g$, we have up to a sign that
    $ \det B  = \pm c \cdot \Res(f,g)  $
    where the extraneous factor $c \in R$ is
    \[ c = \begin{cases}
        f_{m}^{m-n} \ &\mathrm{if} \ m >  n;\\
        g_{n}^{n-m} &\mathrm{if} \ n > m;\\
        1 \ &\mathrm{otherwise}.
    \end{cases}  \]
\end{corollary}
\begin{proof}
    Immediate by Definition \ref{def:res} and Proposition \ref{prop:SFS}.
\end{proof}

We now assume that $R=\FF$ is an algebraically closed field, and prove Lemma \ref{lem:basisforsylvester} and Corollary \ref{cor:basisforbezout} in that context.

\begin{lemma}\label{lem:basisforsylvester}
   Let $f,g,h \in \FF[x]$, with  $h(x)=\gcd(f(x),g(x))$, $m \geq \deg f$ and $n \geq \deg g$. The $(m,n)$-Sylvester matrix of $f$ and $g$ has rank $ m + n - s - \deg \gcd(f,g)$
    where $s := \min(m-\deg f, n - \deg g)$.  Moreover, a basis for $\ker S^\top$ is given by the vectors
    \begin{align*}
         \begin{bmatrix}
            \Gamma_{n}\left(x^i\frac{g(x)}{h(x)}\right)\\
            -\Gamma_{m}\left(x^i\frac{f(x)}{h(x)}\right)
        \end{bmatrix}
    \end{align*}
    for $i = 0,\dots,s+\deg h - 1$, and a basis of $\ker S$ is given by the following vectors:
    \begin{itemize}
        \item The canonical unit vectors $e_i$ for $i = 1,\dots, s$
        \item  For all  $x_0 \in \FF$ such that $f(x_0)=g(x_0)=h(x_0)=0$, the vectors $\Lambda_{m+n}^{(i)}(x_0)$ for $i = 0,\dots,\mu(x_0)-1$ where $\mu(x_0)$ is the multiplicity of $x_0$ as a root of $h(x)$.
    \end{itemize} 
\end{lemma}
Before proceeding with a proof, we note that, among the vectors in the basis of $\ker S$ given in Lemma \ref{lem:basisforsylvester},  the vectors $\{e_i\}_{i=1}^s$ correspond to possible ``shared infinite roots" of $f$ and $g$, whereas the other vectors correspond to common finite roots.
\begin{proof}[Proof of Lemma \ref{lem:basisforsylvester}]
    Let $p(x) \in \FF[x]_{n-1}, q(x) \in \FF[x]_{m-1}$. By item 2 in Proposition \ref{prop:sproperty}, 
     \[ S^\top \begin{bmatrix}
         \Gamma_{n}(p) \\ \Gamma_{m}(q)
     \end{bmatrix} = 0 \iff f(x)p(x)+g(x)q(x)=0 \Leftrightarrow \exists \ a(x) \in \FF[x] : \begin{bmatrix}
            p(x)\\
            q(x)
        \end{bmatrix} = a(x) \begin{bmatrix}
            \frac{f(x)}{h(x)}\\
            -\frac{g(x)}{h(x)}
        \end{bmatrix}. \]
       This induces an isomorphism $\ker S^\top \cong \FF[x]_{s + \deg h -1} \cong \FF^{s+\deg h}$. Note indeed that we look for $p$ and $q$ of degrees at most $m-1$ and $n-1$ respectively, and hence $\deg a \leq s + \deg h - 1$. Thus, $\dim \ker S = s + \deg h$, and the first part of the statement then follows by the rank-nullity theorem.
        By choosing the basis $\{x^i\}_{i=0}^{s+\deg h -1}$ for $\FF[x]_{s + \deg h -1}$, we also obtain the basis of the cokernel $\ker S^\top$ given in the statement.
        
        It remains to verify that the proposed basis of $\ker S$ is indeed a basis.  The given vectors are clearly linearly independent \cite{ait}. Since the $\dim \ker S =m+n-s-\deg h$ by first part of the proof, it suffices to check $Sv=0$ for all vectors $v$ in the proposed basis.  
        If $s > 0$, then the $s$ leftmost columns of $S$ are equal to zero. Indeed, $f_m=f_{m-1}=\dots=f_{m-s+1}=g_n=g_{n-1}=\dots=g_{n-s+1}=0$ because $m-s+1 \geq \deg f +1$ and $n-s+1 \geq \deg g +1$. 
        For the other vectors, by item 1 of Proposition \ref{prop:sproperty} we have
        \[ S \Lambda_{m+n}^{(i)}(x_0)= \left.\begin{bmatrix}
            f(x) \Lambda_n(x)\\
            g(x) \Lambda_m(x)
        \end{bmatrix}^{(i)} \right\vert_{x=x_0} = 0.   \]

\end{proof}

\begin{corollary}\label{cor:basisforbezout}
    Let $B \in \FF^{k \times k}$ be the Bezout matrix of $f$ and $g$, where $k=\max \{ \deg f, \deg g \}$. Then, $B$ has rank $k-\deg \gcd(f,g)$, and a basis for $\ker B$ consists of the vectors $\Lambda_k^{(i)}(x_0)$
where $x_0$ is any shared root of $f, g$ and $0 \leq i < \mu(x_0)$ where $x_0$ is the multiplicity of $x_0$ as a root of
$\gcd(f, g)$.
\end{corollary}
\begin{proof}
Let $S$ be the $(k,k)$-Sylvester matrix of $f$ and $g$. Proposition \ref{prop:SFS} and Frobenius's rank inequality yield
 $ 2 \rank(S) - 2k \leq 2 ~\rank(B),$
implying in turn 
$\dim \ker B \leq \dim \ker S.$ In view of this observation, it suffices to verify $B \Lambda_k^{(i)}(x_0) = 0$ for each allowed choice of $x_0$ and $i$, since (applying Lemma \ref{lem:basisforsylvester} in the case $m=n=k$ and hence $s=0$) these equations
also imply $\dim \ker B \geq \dim \ker S$.

To this goal, suppose that $x_0$ is a root of $\gcd(f, g)$ of multiplicity $\mu(x_0)$ and let $i < \mu(x_0)$. By Proposition \ref{prop:SFS} and Lemma \ref{lem:basisforsylvester}, we have:
\begin{align*}
  0 = S^\top \begin{bmatrix}
      0 & F \\
      -F & 0
  \end{bmatrix} S \Lambda_{2k}^{(i)}(x_0) = \begin{bmatrix}
      0 & B \\
      -B & 0
  \end{bmatrix} \Lambda_{2k}^{(i)}(x_0) = \begin{bmatrix}
      B \Lambda_{k}^{(i)}(x_0)\\
       -B (x^k \Lambda_k)^{(i)}(x_0)
  \end{bmatrix} \Rightarrow B \Lambda_{k}^{(i)}(x_0)=0.
\end{align*}
\end{proof}

    % We may extend the notion to the multivariate case by defining
    % \[ f(x_1,\dots,x_n) = \sum_{i=0}^n f_i(x_2,\dots,x_n) x_1^i \]
    % and same for $g$, and allowing the coefficients to be in $\FF[x_2,\dots,x_n]$. 

\subsection{Smith form and root vectors}
Recall that an $n \times n$ \emph{unimodular} polynomial matrix is a unit of $\FF[y]^{n \times n}$. The classical Theorem \ref{thm:smith} was first proved by H. J. S. Smith \cite{smith} for integer matrices and later by I. Kaplansky \cite{kaplansky} for matrices over any elementary divisor domain (a class of rings that, in particular, includes any principal ideal domain). We note in particular that Theorem \ref{thm:smith}  still holds if $\FF[y]$ is replaced by its localization at a prime $y-y_0$ (see \eqref{eq:localring} below).
\begin{theorem}[Smith's Theorem, \cite{kaplansky,smith}]\label{thm:smith}
    Let $M(y) \in \FF[y]^{m \times n}$. Then, there exist unimodular $U(y) \in \FF[y]^{m \times m}$, $V(y) \in \FF[y]^{n \times n}$ such that $U(y)M(y)V(y)=D(y)$ and $D(y) \in \FF[y]^m$ is diagonal and such that $D_{ii}(y)$ divides $D_{i+1,i+1}(y)$ for all $1 \leq i < \max \{m,n \}$. Moreover, $D(y)$ is uniquely determined by $M(y)$ up to multiplying its diagonal elements by units of $\FF[y]$. The matrix $D(y)$ is called a \emph{Smith form} of $M(y)$, and its nonzero diagonal elements are the \emph{invariant factors} of $M(y)$. A root $y_0 \in \FF$ of an invariant factor of $M(y)$ is called a \emph{finite eigenvalue} of $M(y)$. Moreover, writing a prime factor decomposition $D_{ii}(y)=\prod_{y_j \in \FF} (y-y_j)^{\kappa_{ij}}$ for each invariant factor, the prime factors $y-y_j$ are called \emph{elementary divisors} of $M(y)$ and the positive exponents $\kappa_{ij}$ are called the \emph{partial multiplicities} of the prime $y-y_j$ (or equivalently, of the eigenvalue $y_j$) in $M(y)$.
\end{theorem}

Root vectors, also sometimes called root polynomials, are a useful tool in the theory of polynomial matrices. Definition \ref{def:rootpoly} was given in \cite[Chapter 1]{GLR82} and is valid for the special case of regular polynomial matrices, i.e., square polynomial matrices whose determinant does not identically vanish. For a more general definition see, e.g., \cite{rp1,rp2}.

\begin{definition}[Root vector, \cite{GLR82}]\label{def:rootpoly}
    Let $M(y) \in \FF[y]^{n \times n}$ be regular, i.e., $\det M(y) \neq 0$. A vector $v(y) \in \FF[y]^n$ is called a \emph{root vector of order $\ell \geq 1$} for $M(y)$ at the eigenvalue $y_0 \in \FF$ if
    \begin{itemize}
        \item[(i)] $v(y_0) \neq 0$;
        \item[(ii)] $M(y)v(y)=(y-y_0)^\ell w(y)$ and $w(y_0) \neq 0$.
    \end{itemize}
\end{definition}

We next turn to maximal sets, again specializing the definition in \cite{rp1,rp2} to the regular case.

\begin{definition}[Maximal, complete, and $y_0$-independent sets of root vectors, \cite{rp1,rp2}]
    Let $M(y) \in \FF[y]^{n \times n}$ be regular, fix $y_0 \in \FF$, and suppose that $\{v_i(y) \}_{i=1}^s$ are root vectors for $M(y)$ at $y_0$ of orders $\{ \ell_i\}_{i=1}^s$, with $\ell_1 \geq \dots \geq \ell_s > 0$. We say that $\{v_i(y) \}_{i=1}^s$ is a \emph{$y_0$-independent set} if
   the vectors $\{ v_i(y_0) \}_{i=1}^s \subset \FF^n$ are linearly independent over $\FF$; we say that they are \emph{complete set} if they are a $y_0$-independent set and there does not exist a $y_0$-independent set of root vectors for $M(y)$ at $y_0$ having cardinality $> s$; and we say that they are a \emph{maximal set} if, for all $i=1,\dots,s$, one cannot replace $v_i(y)$ with another root vector $w(y)$ having order $\ell > \ell_i$ and such that $\{ v_1(y),\dots,v_{i-1}(y),w(y),v_{i+1}(y),v_s(y)\}$ is a complete set.
\end{definition}
Complete sets are useful because their cardinality is equal (still assuming that $M(y)$ is regular) to $\dim \ker M(y_0)$ \cite[Theorem 3.10]{rp2}; in turn, this is  equal to the number of partial multiplicities of $y_0$ as an eigenvalue of $M(y)$. Maximal sets are useful because their orders coincide with the partial multiplicities of $y_0$ in $M(y)$. In particular, Theorem \ref{thm:maximalchar} below will be crucial for the goals of this paper.
\begin{theorem}\cite[Theorem 4.1, item 3.]{rp1}\label{thm:maximalchar}
    Let $\{v_i(y) \}_{i=1}^s$ be a complete set of root vectors at $y_0$ for $M(y)$, having orders $\ell_1 \geq \dots \geq \ell_s > 0$. Moreover, suppose that the partial multiplicities of $y_0$ as an eigenvalue of $M(y)$ are $\kappa_1 \geq \dots \geq \kappa_s > 0$. Then, $\ell_i \leq \kappa_i$ for all $i=1,\dots,s$, and the following are equivalent:
    \begin{enumerate}
        \item $\{v_i(y) \}_{i=1}^s$ are a maximal set of root vectors at $y_0$ for $M(y)$;
        \item $\ell_i=\kappa_i$ for all $i=1,\dots,s$;
        \item $\sum_{i=1}^s \ell_i = \sum_{i=1}^s \kappa_i$.
    \end{enumerate}
\end{theorem}

\section{M\"{o}ller indices}\label{sec:indices}

In this section, we define certain indices that describe the shape of every nontrivial closed subspace $V \subset \mathcal{D}$ (see Definition \ref{def:closedsubspace}). We then prove in Theorem \ref{thm:mollerindequiv} that these indices can be easily read off a Gauss basis \cite{3M,3M96} of $V$.

\begin{definition}\label{def:mollerindices}
        Let $\{0\} \neq V \subset \mathcal{D}$ be a nontrivial closed vector subspace, and define
        \[ \beta:=1 + \max \{ i \in \mathbb{N} \ \mathrm{s.t.} \ D_{i0} \in V  \}.  \]
        The positive integers
        \begin{align*}
            \alpha(i):= 1 + \max \{j \in \mathbb{N} \ \mathrm{s.t.} \ \exists \ \phi = D_{ij} + \sum_{\scriptstyle\begin{matrix}
               (k,\ell)\neq(i,j) \\
            \ell \leq j\end{matrix}}\ c_{k\ell} D_{k \ell} \in V \}  , \qquad i=0,\dots,\beta-1
        \end{align*}
        are called the \emph{M\"{o}ller indices with respect to the variable $y$ of $V$}. If $V=V_{(x_0,y_0)}$ is the dual space at $(x_0,y_0) \in \mathcal{V}(I)$ of a finite dimensional ideal $I \subset \FF[x,y]$, the M\"{o}ller indices with respect to $y$ of $V$ are called \emph{M\"{o}ller indices with respect to the variable $y$ of the ideal $I$ at the point $(x_0,y_0)$.} 
    \end{definition}
    
The numbers $\beta \geq 1$ and $\alpha(i) \geq 1$ in Definition \ref{def:mollerindices} are well defined because $V$ is finite-dimensional and nontrivial, and hence the sets $\{i \in \NN \ \mathrm{s.t.} \ D_{i0} \in V\}$ and (for all $0 \leq i < \beta$) $\{j \in \NN \ \mathrm{s.t.} \ \exists \ \phi = D_{ij} + \sum_{
            \ell \leq j, (k,\ell)\neq(i,j)} c_{k\ell} D_{k \ell} \in V\}$ are finite and nonempty; in particular, $D_{00} \in V$. Furthermore, the sequence $\alpha(i)$ is non-increasing in $i$ because  $\phi \in V \Rightarrow A_x \phi \in V$
     \begin{remark}\label{rem:mollerindices}
    For every point $(x_0,y_0)$ in the variety of  the zero-dimensional ideal $\langle f, g \rangle$, the corresponding dual space $V_{(x_0,y_0)}$ depends only on the ideal, and not on the specific choice of the generators $f$ and $g$ \cite{3M,ms}. In turn, Definition \ref{def:mollerindices} only depends on $V_{(x_0,y_0)}$. Therefore, the M\"{o}ller indices are an intrinsic property of a zero-dimensional ideal $\ideal{f,g}$.
\end{remark}
 Visually, if we represent the dual spaces with respect to its leading monomials (in a sense that will be made precise later) as a Young diagram, $\beta$ is the width and $\alpha(i)$ is the height of the $(i-1)$-th column. After developing some required tools, we will illustrate this concept in Figure \ref{fig1}.

    Theorem \ref{thm:beta} below provides another characterization of $\beta$.
\begin{theorem}\label{thm:beta}
    If $V=V_{(x_0,y_0)}$ is the dual space at $(x_0,y_0) \in \mathcal{V}(\ideal{f,g})$, then the number $\beta$ in Definition \ref{def:mollerindices} is equal to the multiplicity of $x_0$ in $\gcd(f(x,y_0),g(x,y_0)) \in \FF[x]$.
\end{theorem}
\begin{proof}
    Let $h(x) := \gcd(f(x,y_0),g(x,y_0))$ and denote by $\mu$ the multiplicity of $x_0$ as a root of $h(x)$. Since $\FF[x]$ is a B\'{e}zout domain, we have $h(x)= a(x) f(x,y_0) + b(x) g(x,y_0)$ for some $a(x),b(x) \in \FF[x] \subset \FF[x,y]$. The fact that $a(x) f(x,y) + b(x) g(x,y) \in \ideal{f,g}$ implies
      \[E_{(x_0,y_0)} D_{i,0} a(x) f(x,y) + b(x)g(x,y) = E_{(x_0,\cdot)} D_{i,0} h(x) = 0  ~\forall i \leq \beta, \]
    and hence, $\mu \leq \beta$.
    
    On the other hand, every element of the ideal can be written as $p(x,y)f(x,y) + q(x,y)g(x,y)$, for some $p,q \in \FF[x,y]$. Thus, using that $\langle f(x,y_0),g(x,y_0)\rangle=\langle h(x) \rangle$,
    \[ E_{(x_0,y_0)} D_{i,0} ~p(x,y)f(x,y) + q(x,y)g(x,y) = E_{(x_0,\cdot)}  D_{i,0} ~p(x,y_0) f(x,y_0) + q(x,y_0)g(x,y_0) \]
  \[  \Rightarrow \exists \ r \in \FF[x] \text{ s.t.} \   E_{(x_0,\cdot)} D_{i,0} ~h(x) r(x) = 0,~ \forall i = 0,\dots,\mu,\]
    showing $\beta \leq \mu$.   
\end{proof}

\begin{remark}\label{rem:swap}
        If in Definition \ref{def:mollerindices} we swap the roles of the first and second indices in the $D_{ij}$, then we can analogously define the M\"{o}ller indices with respect to the variable $x$ of the ideal $I$ at the point $(x_0,y_0)$ or, more generally, of any nontrivial closed subspace $V \subset \mathcal{D}$.
    \end{remark}

    Recall \cite{iva,usingag} that a \emph{monomial order} $<$ is a total order on the set of bivariate monomials $BM:=\{ x^i y^j\}_{(i,j)\in \NN^2}$ such that, for all $u,v,w \in BM$,  $1 \leq u$  and $u < v \implies wu < wv$.
    The bijection $ x^i y^j \mapsto D_{ij}$ shows that every monomial order induces a corresponding order on the $D_{ij}$.  Having fixed a monomial order $<$, the \emph{leading monomial} \cite{iva,usingag} of $0 \neq \phi  = \sum c_{ij} D_{ij} \in \mathcal{D}$ is   $ \LM(\phi) := \max \{ D_{ij} \ \mathrm{s.t.} \ c_{ij} \neq 0 \},$
    where the maximum is with respect to $<$.

    \begin{lemma}\label{lem:LMVisclosed}
        Let $V \subset \mathcal{D}$ be a closed subspace, and consider the set of leading monomials $\LM(V):=\{ \LM(\phi) \mid \phi \in V\}$ with respect to the some monomial order $<$. If $D_{ij} \in \LM(V)$, then $D_{k \ell} \in \LM(V)$ for all $0 \leq k \leq i, 0 \leq \ell \leq j$.
    \end{lemma}
    \begin{proof}
      By assumption, there exists $\phi \in V$ such that $\LM(\phi)=D_{ij}$. Pick any ordered pair $(k,\ell)$ as in the statement, then $\LM(A^{i-k}_x A_y^{j-\ell} \phi)=D_{k\ell} \in \LM(V)$.
    \end{proof}
\begin{definition}
[Gauss basis, \cite{3M,3M96}] \label{G1}
    Let $V \subset \mathcal{D}$ be a closed vector space, and fix a monomial order $<$. If $d:=\dim_\FF V$, a basis $\{G_i\}_{i=1}^d$ of $V$ is called a \emph{ Gauss basis with respect to $<$} if
    \begin{enumerate}
        \item  $\LM(G_i) < \LM(G_j)$ for all $i < j$,
        \item $G_i = \LM(G_i) + \psi_i, \ \psi_i \in \mathcal{D}$, and for all $i,j$, $\LM(G_j)$ is not a monomial in $\psi_i$.
    \end{enumerate}    
\end{definition}

Concretely, a Gauss basis may be constructed by performing Gauss-Jordan elimination on an arbitrary basis \cite{3M}. Among other reasons, Gauss bases were introduced in \cite{3M} because they suffice for the purpose of studying $\LM(V)$; we record this property in Proposition \ref{prop:wasitemi}.
\begin{prop}\label{prop:wasitemi}
Let $V \subset \mathcal{D}$ be a closed vectors subspace with dimension $d$, fix a monomial order $<$, let $\LM(V)$ be as in Lemma \ref{lem:LMVisclosed} and let $\{G_i\}_{i=1}^d$ be a Gauss basis of $V$ with respect to $<$. Then, $\LM(V)=\{ \LM(G_i)  \}_{i=1}^d$.
\end{prop}
\begin{proof}
    For every $0\neq \phi \in V$, there exist $c_1,\dots,c_d \in \FF$, not all zero and such that
    $ \phi = \sum_{i=1}^d c_i G_i \Rightarrow \LM(\phi)=\LM(G_s)$
    where, by property 1. in Definition \ref{G1}, $ s= \max \{ i : c_i \neq 0 \}$. This shows $\LM(V) \subseteq \{\LM(G_i)\}_{i=1}^d$; the reverse inclusion is obvious since $\{G_i\}_{i=1}^d \subset V$. 
\end{proof}

We now restrict our attention to the \emph{lex order with $x < y$} \cite{iva,usingag}:
\[ D_{ij} < D_{k \ell} \iff (j < \ell) \text{ or } (j= \ell) \wedge (i < k). \]

Theorem \ref{thm:mollerindequiv} below shows that the M\"{o}ller indices with respect to $y$ of $V$ can be easily computed from a Gauss basis of $V$ with respect to the lex order with $x < y$. We note in passing that a Gauss basis with respect to the lex order with $y<x$ would instead provide the M\"{o}ller indices with respect to $x$.

\begin{theorem} \label{thm:mollerindequiv}
    Let $\{0\} \neq V \subset \mathcal{D}$ be a nontrivial closed vector subspace with dimension $d=\dim_\FF V>0$, let $\{ G_i \}_{i=1}^d$ be a Gauss basis of $V$ with respect to the lex order with $x<y$, let the M\"{o}ller indices with respect to $y$ of $V$ be $\alpha(0) \geq \dots \geq \alpha(\beta-1)$, and let $\LM(V)$ be as in Lemma \ref{lem:LMVisclosed}. Then, 
    \begin{itemize}
        \item[(i)]  $G_k= D_{k-1,0}$ for all $1 \leq k \leq \beta$; moreover, if $k > \beta$ and $D_{ij}=\LM(G_k)$ then $j>0$;
        \item[(ii)]  If we partition
    \[  \LM(V) = \bigcup_{i=0}^{\gamma-1} \mathcal{L}_i, \qquad \mathcal{L}_i : =   \{ D_{ij} \in \LM(V) \ \mathrm{for} \ \mathrm{some} \ j \in \NN  \},\]
    then $\gamma=\beta \geq 1$ and the cardinality of the set $\mathcal{L}_i$ is $\# \mathcal{L}_i = \alpha(i) \geq 1$. 
    \end{itemize}
\end{theorem}
\begin{proof}
\begin{itemize}
    \item[(i)] By Proposition \ref{prop:wasitemi} and Definition \ref{G1}, for all $0 \leq k < \beta$, $D_{k0} \in \{ G_i\}_{i=1}^d$, and therefore $G_{k+1}=D_{k0}$.  Suppose now  that $D_{i0}=\LM(G_j)$ for $j>\beta$, then by Lemma \ref{lem:LMVisclosed} and Proposition \ref{prop:wasitemi} $D_{\beta0} \in \LM(V) = \{\LM(G_i)\}_{i=1}^d$. Again by Definition \ref{G1}, this in turn yields $D_{\beta0}=G_{\beta+1} \in V$, contradicting the definition of $\beta$.
   \item[(ii)]  Item (i) above implies $D_{i0} \in V \Leftrightarrow D_{i0} \in \LM(V)$, and hence $\beta=\gamma$.

    Fix now $0 \leq i < \beta$. By Definition \ref{def:mollerindices}, there exists $\phi \in V$ of the form
    \[ \phi = D_{i,\alpha(i)-1} + \sum_{\scriptstyle\begin{matrix}
        (k,\ell)\neq(i,j) \\
        \ell < \alpha(i)
    \end{matrix}} c_{k \ell} D_{k \ell} \Rightarrow \LM(\phi)=D_{a,\alpha(i)-1}  \ \mathrm{for} \ \mathrm{some} \ a \geq i. \]
 Lemma \ref{lem:LMVisclosed} implies that $\# \mathcal{L}_i$ is a non-increasing function of $i$. Therefore, 
     $\alpha(i) \leq \# \mathcal{L}_a \leq \# \mathcal{L}_i.$ Conversely, by Propsition \ref{prop:wasitemi}, $D_{i,\# \mathcal{L}_i-1}=\LM(G_b)$ for some $1 \leq b \leq d$. Therefore, 
     \begin{equation}\label{eq:katsominua}
      \ G_b = D_{i, \# \mathcal{L}_i-1} + \sum_{k < i} c_k D_{k, \# \mathcal{L}_i-1} + \sum_{\scriptstyle\begin{matrix}
        k \in \NN \\
        \ell < \# \mathcal{L}_i-1
    \end{matrix}} c_{k \ell} D_{k \ell} \in V,      
     \end{equation}  
    and hence $\# \mathcal{L}_i \leq \alpha(i)$.
\end{itemize}

\end{proof}

\begin{corollary}\label{cor:etna}
     Let $\ideal{f,g}$ be zero-dimensional, let $P=(x_0,y_0) \in \mathcal{V}(\ideal{f,g})$, and let $V_{P}$ be the dual space of $\ideal{f,g}$ at $P$. Denoting by  $\alpha(0) \geq \dots \geq \alpha(\beta-1)$ the associated M\"{o}ller indices with respect to $y$, then $\dim V_{P}=\sum_{i=0}^{\beta-1} \alpha(i)$.
\end{corollary} 
Item (ii) in Theorem \ref{thm:mollerindequiv}, and in particular \eqref{eq:katsominua}, prompt Definition \ref{def:lv}.
\begin{definition}[Leading vector] \label{def:lv} Let $\mathcal{G}=\{G_i\}_{i=1}^d$ be a Gauss basis with respect to the lex order with $x<y$ of the nontrivial vector subspace $V \subset \mathcal{D}$, and let  the M\"{o}ller indices of $V$ with respect to $y$ be $\alpha(0) \geq \dots \geq \alpha(\beta-1)$. For $i=0,\dots,\beta-1$, the \emph{$i$-th leading vector of $\mathcal{G}$} is $\phi_i : = G_{b_i}$, where $b_i$ is such that $\LM(G_{b_i})=D_{i,\alpha(i)-1}$.    
\end{definition}

\begin{remark}
We can refine the expression \eqref{eq:katsominua} for the leading vectors of $\mathcal{G}$ by observing that, for every $0 \leq \ell \leq i$, $ \LM(A_x^\ell G_{b_i}) = D_ {i-\ell,\alpha(i)-1} \in \LM(V)=\LM(\mathcal{G})$.
Hence, by Property 2. in Definition \ref{G1}, the leading vector $\phi_i$ must have the form
\begin{equation}\label{eq:mirame}
    \phi_i=G_{b_i} = D_{i,\alpha(i)-1} + \sum_{\scriptstyle\begin{matrix}
        k \in \NN \\
        \ell < \# \mathcal{L}_i-1
    \end{matrix}} c_{k \ell} D_{k \ell}.
\end{equation}
\end{remark}

\begin{example}\label{ex1}
Let $\FF=\mathbb{C}$, $f=(x+y)^2$ and $g=x^3-y^3$. Then $\mathcal{V}(\ideal{f,g})=\{(0,0)\}$, and a Gauss basis with respect to the lex order with $x<y$ of $V_{(0,0)}$ is $\{G_i\}_{i=1}^6$, where $G_1=D_{00}, G_2=D_{10}, G_3=D_{01},$ \[ G_4=D_{11}-2 D_{20}, G_5=D_{02}-D_{20},G_6=D_{03}-\frac13 D_{12} - \frac{1}{3} D_{21} + D_{30}.\]
Hence, the associated M\"{o}ller indices with respect to $y$ are $4$ and $2$, because
\[ \mathcal{L}_0=\{ G_1,G_3,\LM(G_5),\LM(G_6) \}, \qquad \mathcal{L}_1 = \{ G_2, \LM(G_4) \}.   \]
The leading vectors of this Gauss basis are $\phi_0=G_6=D_{03}-\frac{1}{3}D_{12} - \frac{1}{3} D_{21} + D_{30}$ and $\phi_1=G_4=D_{11}-2D_{20}$.
\begin{figure}[h!]
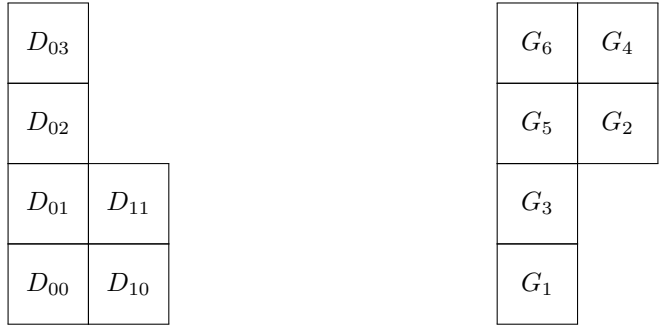

    \centering
    \ytableausetup{boxsize=3em}
       \begin{ytableau}
    D_{03}  \\
    D_{02}  \\
    D_{01} & D_{11} \\
    D_{00} & D_{10}
\end{ytableau}
\hspace{4cm}
\ytableausetup{boxsize=3em}
   \begin{ytableau}
    G_6 & G_4 \\
    G_5 & G_2 \\
    G_3 \\
    G_1
\end{ytableau}
    \caption{Visualization of Example \ref{ex1} by Young diagrams with width $\beta$. \emph{Left}: Leading monomials, French style diagram. The $i$-th column corresponds to the subset $\mathcal{L}_{i-1}$ in the partition of $\LM(V)$, and thus its height is $\alpha(i)$. Moving leftwards (resp. downwards) correspond to the action of $A_x$ (resp. $A_y$).  \emph{Right}: Corresponding vectors in the Gauss basis $\{ G_i\}_{i=1}^6$, English style diagram. The height of each column still captures the M\"{o}ller indices with respect to $y$. The leading vectors ($\phi_0=G_6$, $\phi_1=G_4$) appear on top of each column.}
    \label{fig1}
\end{figure}
\end{example}

\section{Smith form of Sylvester matrices}\label{sec:sylvestersmith}

Let $\langle f,g \rangle$ be a zero-dimensional ideal of $\FF[x,y]$. Our aim is to describe the Smith form of the Sylvester matrix of $f$ and $g$, seen as elements of $\FF[x][y]$, in terms of the M\"{o}ller indices with respect to $y$ of $\langle f,g \rangle$.

\subsection{Case of coprime leading coefficients}

It is convenient to first assume that the leading coefficients of $f$ and $g$ (as polynomials in $x$) do not have non-unital common divisors in $\FF[y]$. This ensures that there are no infinite common roots; the main result of this subsection is Theorem \ref{thm:mainresultnoinf}, that describes all the invariant factors of $S(y)$ in this case. Later, in Theorem \ref{thm:mainresult}, we will generalize Theorem \ref{thm:mainresultnoinf} by removing the assumption that the leading coefficients are coprime. We begin with the auxiliary Lemma \ref{lem:bivariatedivision}.
    
\begin{lemma}\label{lem:bivariatedivision}
    Let $R$ be an integral domain and pick $p,q \in R[x]$ such that the leading coefficient of $p$ is a unit of $R$. Let $r \in \langle p,q \rangle$ and suppose that $\deg r < \deg p + \deg q$. Then, there exist $u,v \in R[x]$ such that (i) $\deg u < \deg q$ (ii) $\deg v < \deg p$ (iii) $r=up+vq$.
\end{lemma}
\begin{proof}
    Let $\tilde{u}, \tilde{v} \in R[x]$ be such that $r=\tilde{u}p+\tilde{v}q$; their existence is guaranteed because $r \in \langle p,q \rangle$. By monic polynomial division \cite[pp. 128-129]{jacobson}, there exist $u,t \in R[x]$ such that $\tilde{v}=tp+v$ and $\deg v < \deg p$. Hence, $r=(\tilde{u}+tq)p+vq$. Define $u:=\tilde{u}+tq$ and suppose for a contradiction that $\deg u \geq \deg q$. Then, 
\[ \deg v + \deg q <  \deg p + \deg u \Rightarrow \deg r =\deg u + \deg p \geq \deg q + \deg p,   \]
    violating the assumptions.
\end{proof}

The next step is to apply the structure theorem and observe in Theorem \ref{thm:wehatealggeom} that, when there are no shared infinite roots, the modular decomposition of the Sylvester matrix describes the algebraic structure of the quotient $\FF[x,y]/\langle f,g \rangle$ as a module over the ring $\FF[y]$. We first do this locally at each eigenvalue, in Theorem \ref{thm:welovealggeom}.

\begin{theorem}\label{thm:welovealggeom}
    Let $f,g \in \FF[x,y]$, and suppose that
    \begin{enumerate}
        \item[(i)] the ideal $I := \ideal{f,g}$ is zero-dimensional
        \item[(ii)] the leading coefficients in $x$ of $f$ and $g$ are coprime in $\FF[y]$.
    \end{enumerate}
     Denote the localization of $\FF[y]$ at the prime ideal $y-y_0$ by
    \begin{equation}\label{eq:localring}
        R := \FF[y]_{(y-y_0)} = \{ a(y)/b(y) \mid a,b \text{ coprime}, \  b(y_0)\neq 0\},
    \end{equation} 
    and the Sylvester matrix associated with $f$ and $g$ by $S(y)$. Suppose that the local Smith form of $S(y)$ at $y_0$ (that is, the Smith form of $S(y)$ over the ring $R)$ has the nontrivial invariant factors $(y-y_0)^{\kappa_i}$, $i=1\,\dots,r$. Then, as $R$-modules, $R[x]/I \cong \oplus_{i=1}^r R/\ideal{(y-y_0)^{\kappa_i}}$.
\end{theorem}
\begin{proof}
    Let $\mathcal{S}$ be the Sylvester map defined as in Proposition \ref{prop:sproperty}, seen as a linear map between $R$-modules and choosing $m=\deg_x f$, $n=\deg_x g$. 
    Since $\mathcal{S}$ is represented, in the monomial basis, by $S^\top(y)$, we have by the structure theorem for finitely generated modules 
    \[ \mathrm{coker} \  \mathcal{S} := R[x]_{m+n-1}/\mathrm{Im} \ \mathcal{S}  \cong \bigoplus_{i=1}^r R/\langle (y-y_0)^{k_i}\rangle.   \]
    What remains to be proven is that $\mathrm{coker} \  \mathcal{S}\cong R[x]/I$ as $R$-modules. Set for simplicity of notation $M:=\mathrm{coker} \  \mathcal{S}$ and $N:=R[x]/\langle f,g \rangle$, and consider the map $\phi$ defined, for all $a \in R[x]_{m+n-1}$, as
\[ \phi:  M \to N, \qquad \phi([a]_M)=[a]_N. \]
It is obvious that $\phi$ is $R$-linear; we must show bijectivity. 
Since the leading coefficients of $f,g$ are coprime over $\FF[y]$ (and hence, a fortiori, over $R$), at least one of them is a unit of $R$; suppose without loss of generality that this is true for the leading coefficient of $f$.

Let $a \in R[x]_{m+n-1}$, be such  that $\phi(a) = [0]_N$, that is, $a \in I$. 
Since the leading coefficient of $f$ is a unit of $R$, by Lemma \ref{lem:bivariatedivision} there exist $u \in R[x]_{n-1}$ and $v \in R[x]_{m-1}$ such that $a=uf+vg$. 
Hence, $a \in \mathrm{Im} \ \mathcal{S}$ and $[a]_M = [0]_M$. 

Now take an arbitrary $b \in R[x]$. By polynomial division we can
    uniquely write $b=c f + d$ with $\deg_x d < m \leq m+n-1$, that is, $d \in R[x]_{m+n-1}$. Hence, $[d]_M$ is well defined and $[b]_N = [d]_N = \phi([d]_M)$. This shows that $\phi$ is surjective.
\end{proof}

% \begin{remark}   
%     If for a fixed $y_0 \in \FF$ there is unique $x_0 \in \FF$ such that $(x_0,y_0) \in \mathcal{V}$, then Theorem \ref{thm:welovealggeom} implies $V_{(x_0,y_0)} \cong \oplus_{i=1}^r R/\ideal{(y-y_0)^{\kappa_i}}$, that is $\dim V_{(x_0,y_0)} = \sum_{i=1}^r \kappa_i$.
%     Otherwise, we have $\oplus_{x_0\in \FF} V_{(x_0,y_0)} \cong \oplus_{i=1}^r R/\ideal{(y-y_0)^{\kappa_i}}$, and 
%     \[ \sum_{x_0 \in \FF} \dim V_{(x_0,y_0)} = \sum_{x_0 \in \FF} \sum_{i=0}^{\beta^{\{x_0\}}-1} \alpha^{\{x_0\}}(i)  = \sum_{i=1}^r  \kappa_i.  \]
% \end{remark}

\begin{theorem}\label{thm:wehatealggeom}
    Under the same assumptions as Theorem \ref{thm:welovealggeom}, let $S(y)$ have (over $\FF[y]$) the nontrivial invariant factors $s_1(y),\dots,s_r(y)$. Then, as $\FF[y]$-modules, $\FF[x,y]/\ideal{f,g} \cong \oplus_{i=1}^r \FF[y]/\ideal{s_i(y)}$.
\end{theorem}
\begin{proof}
    It follows immediately from Theorem \ref{thm:welovealggeom}. Indeed, since $\FF[y]$ is a PID, the global invariant factors of $S(y)$ are the product of the local factors at each prime ideal $y-y_0$. 
\end{proof}

To state Corollary \ref{cor43}, recall that, given a prime $y-y_0 \in \FF[x]$ and an element $q \in \FF[x]$, the local $p$-adic valuation $\nu_{y_0}$ is $\nu_{y_0}(q)=e$ if $(y-y_0)^e$ divides $q$ but $(y-y_0)^{e+1}$ does not.

\begin{corollary}\label{cor43}
    Under the assumptions of Theorems \ref{thm:welovealggeom}, 
    \[ \deg_y \det S(y) = \deg_y \Res_x(f,g) = \dim_{\FF} \FF[x,y]/\ideal{f,g}. \]
    Morever, for all $y_0 \in \FF$ and if $R$ is the local ring \eqref{eq:localring}, 
     \[ \nu_{y_0} (\det S(y) )= \nu_{y_0} (\Res_x(f,g)) = \dim_{\FF} R/\ideal{f,g}. \]
\end{corollary}
\begin{proof}
   We prove the first statement, as the argument for the second is analogous. 
   
   $\FF[y]/\langle s_i(y) \rangle \cong \FF^{\deg_y s_i(y)}$
    as $\FF$-vector spaces. Let $N:=\sum_{i=0}^g \deg_y s_i(y) = \deg_y \det S(y)=\deg_y \mathrm{Res}(f,g). $ Theorem \ref{thm:wehatealggeom} then yields the vector space isomorphism $\FF[x,y]/\langle f,g \rangle \cong \FF^N$.
\end{proof}

\begin{corollary} \label{cor:resdeg}
  Consider the zero-dimensional ideal $I:=\langle f,g \rangle$. For every point $P \in \mathcal{V}(I) \subset \FF^2$, let $V_{P}$ denote the dual space of $I$ at the point $P$, and let $\alpha_P(0),\dots,\alpha_P(\beta_P-1)$ be the M\"{o}ller indices with respect to $y$ of $I$ at $P$. Then, under the same assumptions of Theorem \ref{thm:welovealggeom},
   \[ \deg_y \det S(y) = \deg_y \Res_x(f,g) = \sum_{P \in \mathcal{V}(I)} \dim V_{P} = \sum_{P \in \mathcal{V}(I)} \sum_{i=0}^{\beta_P-1} \alpha_P(i). \]
   Moreover, for all $y_0 \in \FF$, let $\mathcal{W}_{y_0}$ be the subset of $\mathcal{V}(I)$ containing all and only the points of the form $(x_i,y_0)$ for some $x_i \in \FF$. Then,
 \[ \nu_{y_0}(\det S(y)) = \nu_{y_0}(\Res_x(f,g)) = \sum_{P \in \mathcal{W}_{y_0}} \dim V_{P} = \sum_{P \in \mathcal{W}_{y_0}} \sum_{i=0}^{\beta_P-1} \alpha_P(i). \]
\end{corollary}
\begin{proof}
The statement follows by Corollary \ref{cor:etna}, Corollary \ref{cor43}, \cite[Theorem 1]{ms} and \cite[Theorem 2.6]{3M}.
\end{proof}

To simplify the notation in what follows, we now define certain linear operators.

\begin{definition}
    Given $y_0 \in \FF$ and $\alpha \in \mathbb{N}$, the linear operator $R_{y_0,\alpha} : \mathcal{D} \to \FF[x,y]^{m+n}$ is defined by
    \[ R_{y_0,\alpha} D_{ij} := \begin{cases}
        (y-y_0)^{\alpha - j}\Lambda^{(i)}_{m+n}(x) & \mathrm{if} \ \alpha \geq j;\\
        0&\text{otherwise}.
    \end{cases} \]
\end{definition}

Note that, for all $j,k \leq \alpha$ and denoting by $S(y)$ the Sylvester matrix of $f$ and $g$, we have
    \begin{align}
        &E_{(x_0,y_0)} D_{0,k} ~S(y) R_{y_0,\alpha} D_{ij}  \notag \\
        &= E_{(x_0,y_0)} D_{0,k}(y-y_0)^{\alpha - j} S(y) \Lambda^{(i)}_{m+n}(x) \notag \\
        &= \sum_{\ell = 0}^k \underbrace{E_{(x_0,y_0)} D_{0,k-\ell} S(y) \Lambda^{(i)}_{m+n}(x) D_{0,\ell} (y-y_0)^{\alpha - j}}_{=0 \ \mathrm{unless}  \ \ell=\alpha-j} \notag \\
        &=  E_{(x_0,y_0)} A_y^{\alpha-k} D_{i,j} \begin{bmatrix}
            f(x,y) \Lambda_{n}(x) \\
            g(x,y) \Lambda_{m}(x)
        \end{bmatrix}.\label{eq:simplifyRk}
    \end{align}
    In particular, when $k = \alpha$,
    \begin{align} \label{eq:simplifyR}
        E_{(x_0,y_0)} D_{0,\alpha} ~S(y) R_{y_0,\alpha} D_{ij} =  E_{(x_0,y_0)} D_{i,j} \begin{bmatrix}
            f(x,y) \Lambda_{n}(x) \\
            g(x,y) \Lambda_{m}(x)
        \end{bmatrix}.
    \end{align} 

In Lemma \ref{lem:rootpolys}, we construct root vectors for $S(y)$ at the eigenvalue $y_0$ by means of the operators $R_{y_0, \alpha}$.
    
\begin{lemma}\label{lem:rootpolys}
Let $\langle f ,g \rangle$ be zero-dimensional and fix $y_0 \in \FF$. For every $x_0 \in \FF$ such that $(x_0,y_0) \in \mathcal{V}(\langle f,g \rangle)$,  denote by $\phi_i$ the $i$-th leading vector of any Gauss basis, with respect to the lex order with $x<y$, of the dual space of $\langle f,g \rangle$ at $(x_0,y_0)$, and by $\alpha(0),\dots,\alpha(\beta-1)$ the M\"{o}ller indices with respect to $y$ of $\langle f,g \rangle$ at $(x_0,y_0)$. Then, for all $i=0,\dots,\beta-1$ 
    \[ r_i^{(x_0)}(y) := E_{(x_0,\cdot)} R_{y_0,\alpha(i)-1}(\phi_i) \]
    is a root vector of $S(y)$ at $y_0$ and it has order $\geq \alpha(i)$.
\end{lemma}
\begin{proof}
We check that $r_i(y)^{(x_0)}$ satisfies the two items in Definition \ref{def:rootpoly}, with $\ell \geq \alpha(i)$.
\begin{itemize}
    \item[(i)] Upon inspection of \eqref{eq:mirame}, we see that the only term in $\phi_i$ whose second index is $\geq \alpha(i)-1$ is $\LM(\phi_i)=D_{i,\alpha(i)-1}$, with coefficient $1$. Hence, $r_i^{(x_0)}(y_0) = \Lambda_{m+n}^{(i)}(x_0) \neq 0$. 
    \item[(ii)]  Our strategy is to apply Proposition \ref{prop:primechar}. By \eqref{eq:simplifyR} we have
    \begin{align*}
        E_{(\cdot ,y_0)} D_{0,\alpha(i)-1} S(y) r_i^{(x_0)}(y) = %E_{x_0,y_0}  D_{0,\alpha(i)} S(y) R_{y_0,\alpha(i)}(\phi_i) =
        E_{x_0,y_0} \phi_i \begin{bmatrix}
            f(x,y) \Lambda_{n}(x) \\
            g(x,y) \Lambda_{m}(x)
        \end{bmatrix} = 0,
    \end{align*}
    since $\phi_i \in V_{(x_0,y_0)}$.
    In other words, the  $(\alpha(i)-1)$-th Hasse derivative with respect to $y$ of $S(y) r_i^{(x_0)}(y)$ vanishes at $y_0$.
    Furthermore, the fact that the dual space is a closed subspace implies that the Hasse derivatives of lower order also all vanish, because $A_y^k \phi_i \in V_{(x_0,y_0)}$ for all $k=0,\dots,\alpha(i)-1$, and hence \eqref{eq:simplifyRk} implies that
    \begin{align*}
        E_{(\cdot,y_0)} D_{0,k} S(y) r_i^{(x_0)}(y) = E_{(x_0,y_0)} A_y^{\alpha-k} \phi_i \begin{bmatrix}
            f(x,y) \Lambda_{n}(x) \\
            g(x,y) \Lambda_{m}(x)
        \end{bmatrix}.
    \end{align*}

\end{itemize}
\end{proof}

\begin{lemma}\label{lem:wew}
    Let $r^{(x_0)}_i(y)$ be defined as in Lemma \ref{lem:rootpolys}. Let \begin{equation}\label{eq:maximalset}
       \mathcal{R}: = \bigcup_{x_0 : (x_0,y_0) \in \mathcal{V}(\langle f,g\rangle)} \{r_i^{(x_0)}(y)\}_{i=0}^{\beta(x_0)-1},
    \end{equation}
    where the outer union is taken over all common roots $x_0$ of $f(x,y_0)$ and $g(x,y_0)$ and, for each such $x_0$, $\beta(x_0)$ denotes the positive integer $\beta$ associated with the dual space $V_{(x_0,y_0)}$ of $\ideal{f,g}$. Then, $\mathcal{R}$ is a complete set of root polynomials for $S(y)$ at $y_0$.
\end{lemma}
\begin{proof}
    The set $\mathcal{R}$ is $y_0$-independent, because the vectors $r_i^{x_0}(y_0) =  \Lambda_{m+n}^{(i)}(x_0)$ are the columns of a confluent Vandermonde matrix, which has full rank \cite{ait}.

    By Theorem \ref{thm:beta}, $\beta(x_0)$ is the multiplicity of $x_0$ as a root of $\gcd(f(x,y_0),g(x,y_0))$.
    Hence, the cardinality of $\mathcal{R}$ is equal to $\sum_{x_0} \beta(x_0)$.
    In turn, this coincides with the dimension of $\ker S(y_0)$, which is equal to $\deg \gcd(f(x,y_0),g(x,y_0))$ by Lemma \ref{lem:basisforsylvester}. Hence, by \cite[Theorem 3.10]{rp2}, $\mathcal{R}$ is a complete set.
\end{proof}

\begin{lemma}\label{lem:wow}
    The set $\mathcal{R}$ defined in \eqref{eq:maximalset} is a maximal set of root vectors at $y_0$ for $S(y)$.
\end{lemma}
\begin{proof}

    Denote by $\ell_{i}^{\{x_k\}}$ the orders of the root polynomial $r_i^{(x_k)}$ in \eqref{eq:maximalset}, by $\kappa_j$ the partial multiplicities in $S(y)$ of the eigenvalue $y_0$.
    Defining $\mathcal{W}_{y_0}$ as in Corollary \ref{cor:resdeg}, by Theorem \ref{thm:maximalchar} and the local statement in Corollary \ref{cor:resdeg},
    \[ \nu_{y_0} (\det S(y)) = \sum_{j} \kappa_j   = \sum_{x_k \in W_{y_0}} \sum_{i=0}^{\beta(x_k) -1} \alpha_{x_k}(i) \leq \sum_{x_k \in W_{y_0}} \ell_{i}^{\{x_k\}} \leq \sum_i \kappa_i. \]
    Theorem \ref{thm:maximalchar} then implies maximality of the set.
\end{proof}

Lemma \ref{lem:wow} immediately implies Theorem \ref{thm:mainresultnoinf}, which connects the partial multiplicities of the Sylvester matrix $S(y)$ associated with $f$ and $g$ and the M\"{o}ller indices of the ideal $\langle f, g \rangle$.

\begin{theorem}\label{thm:mainresultnoinf}
    Let $S(y) \in \FF[y]^{(m+n) \times (m+n)}$ be the Sylvester matrix of $f,g \in \FF[x,y]$, where $m=\deg f$ and $n=\deg g$. Suppose that 
    \begin{enumerate}
        \item[(i)] $ \ideal{f,g}$ is zero-dimensional and
        \item[(ii)] the leading coefficients in $x$ of $f$ and $g$ are coprime over $\FF[y]$.
    \end{enumerate}
    Fix $y_0 \in \FF$, let $x_i \in \FF$, $i=1,\dots,\ell$, denote all the elements of $\FF$ such that $f(x_i,y_0)=g(x_i,y_0)=0$, and let $\alpha_i(0),\dots,\alpha_i(\beta_i-1)$ be the M\"{o}ller indices with respect to $y$ of the the ideal $\langle f,g \rangle$ at the point $P_i=(x_i,y_0)$. Then, the partial multiplicities of $y_0$ in $S(y)$ are
\[ \alpha_{i}(k), \quad k=0,\dots,\beta_{i}-1, \quad i=1,\dots,\ell.  \]
In particular, the geometric multiplicity of $y_0$ is equal to
\[ \sum_{i=1}^\ell \beta_{i}  \]
and the algebraic multiplicity of $y_0$ is equal to
    \[  \sum_{i=1}^\ell  \sum_{k=0}^{\beta_{i}-1} \alpha_{i}(k) .\]
\end{theorem}

\begin{example}
 Let $\FF=\mathbb{C}$ and let $f,g$ be as in Example \ref{ex1}. Recall that we computed that their only intersection is $(0,0)$, for which $\beta=2$ and $\alpha(0)=4 > \alpha(1)=2$. The Sylvester matrix of $f$ and $g$ is
 \[ S(y) = \begin{bmatrix}
     1 & 2 y & y^2 & 0 & 0\\
     0 & 1 & 2y & y^2 & 0\\
     0 & 0 & 1 & 2y & y^2 \\
     1 & 0  & 0 & -y^3 & 0\\
     0 & 1 & 0 & 0 & -y^3
 \end{bmatrix}.\]
 The $3 \times 3$ leading principal minor of $S(y)$ is equal to $1$, implying that the first three invariant factors of $S(y)$ are $1$. Call now $\sigma_4(y),\sigma_5(y)$ the fourth and fifth invariant factors of $S(y)$, respectively. Denoting the GCD over $\mathbb{C}[y]$ of the entries of the polynomial matrix $M(y)$ by $\gcd(M(y))$, and denoting by the symbol $\dot{=}$ equalities that hold up to multiplication by a unit of $\mathbb{C}[y]$ (that is, a nonzero complex number), we notice that \[ \sigma_4(y)\sigma_5(y) \dot{=}\det S(y)\dot{=}y^6 \quad \mathrm{and} \quad \sigma_4(y)\dot{=}\gcd(\mathrm{adj}(S(y)))\dot{=}y^2.\] Hence, the Smith form of $S(y)$ is indeed $I_3 \oplus y^2 \oplus y^4$, as predicted by Theorem \ref{thm:mainresultnoinf}. Moreover, the maximal set \eqref{eq:maximalset} of root polynomials of $S(y)$ at $0$ is (denoting by $\{e_j\}_{j=1}^5$ the canonical basis of $\mathbb{C}^5$)
 $r_0(y)=e_5 - \frac{y}{3} e_4 - \frac{y^2}{3} e_3 + y^3 e_2$, $r_1(y)=e_4 - 2 y e_3$. Indeed, \[ S(y)r_0(y)=y^4 \frac{5e_1 +e_4}{3} , \qquad S(y)r_1(y)=y^2(-3e_2 - 2 y e_1 - y e_4).   \]
\end{example}

\subsection{Infinite intersections of bivariate polynomials}

When $f$ and $g$ do not satisfy the assumption (ii) in Theorem \ref{thm:welovealggeom}, then neither Theorems  \ref{thm:welovealggeom} and \ref{thm:wehatealggeom} nor, a fortiori, Corollaries \ref{cor43} and \ref{cor:resdeg} hold any longer. For example, consider $\FF=\mathbb{C}$, $f=xy+1$ and $g=x y^2-1$. It is easy to verify that the Sylvester matrix of $f$ and $g$ is $S(y)=\begin{bmatrix}
    y & 1\\
    y^2 & -1
\end{bmatrix}$, whereas a reduced Gr\"{o}bner basis of $\langle f, g \rangle$ is $\{x-1,y+1\}$. Thus, $\dim_\FF \FF[x,y]/\langle f,g \rangle = 1 \neq 2 = \deg \Res_x(f,g)$. The reason of this discrepancy is that, in this example, $f$ and $g$ share a ``common root at $(x,y)=(\infty,0)$", in the sense that their degrees in $x$ both drop when $y=0$.

We now formalize this idea by extending the definitions of roots, intersection multiplicity, and M\"{o}ller indices, with the ultimate goal of going around this difficulty.

\begin{definition}
Let $f \in \FF[y][x]$ have degree $m$ in $x$. Writing
$f(x,y) = \sum_{i=0}^{m} f_i(y) x^i,$ the \emph{reversal} with respect to $x$ of $f$ is
\[ \mathrm{rev}_x f(x,y) = \sum_{i=0}^m f_{m-i}(y) x^i. \]
\end{definition}

\begin{definition}\label{def:projectiveintersection}
    Let $f,g \in \FF[x,y]$ and fix $y_0 \in \FF$. Then:
    \begin{enumerate}
        \item We say that $(\infty,y_0)$ is an infinite intersection of $f$ and $g$, or equivalently a common root of $f$ and $g$ at infinity, if $(0,y_0) \in  \mathcal{V}(\langle \mathrm{rev}_x f, \mathrm{rev}_x g\rangle)$;
        \item The intersection multiplicity of an infinite intersection $(\infty,y_0)$ of $f$ and $g$ is the intersection multiplicity of $(0,y_0)$ as a point in $\mathcal{V}(\langle \mathrm{rev}_xf, \mathrm{rev}_x g \rangle)$;
        \item The M\"{o}ller indices with respect to $y$ of $f$ and $g$ at $(\infty,y_0)$ are the M\"{o}ller indices with respect to $y$ of $\langle \mathrm{rev}_xf, \mathrm{rev}_x g \rangle$ at $(0,y_0)$.
    \end{enumerate}
\end{definition}

\begin{remark}
    Infinite intersections of $f$ and $g$ having the form $(x_0,\infty)$ can be defined analogously by considering instead the reversals of $f$ and $g$ with respect to $y$.
\end{remark}

Note that $(\infty,y_0)$ is a common root at infinity of $f$ and $g$ if and only if the leading coefficients in $x$ of $f$ and $g$ have a common root at $y_0$. The main advantage of Definition \ref{def:projectiveintersection} is that it allows us to obtain analogues of Theorem \ref{thm:wehatealggeom}, Corollary \ref{cor43} and Corollary \ref{cor:resdeg}, and  thus to generalize Theorem \ref{thm:mainresultnoinf} to the case of non-coprime leading coefficients. On the other hand, unlike what happens at any finite intersection point in $\FF^2$ and justifying our lexical choice in item 3 of Definition \ref{def:projectiveintersection}, the M\"{o}ller indices of $f$ and $g$ at an infinite intersection are \emph{not} uniquely determined by $\langle f,g \rangle$, but depend on the choice of the generators $f$ and $g$. To see this, consider again the example $\FF=\mathbb{C}$, $f=xy+1, g=xy^2-1, p=x-1, q=y+1$. Then, $\langle f,g \rangle = \langle p, q \rangle$, but \[ \langle x+y,x-y^2 \rangle = \langle \mathrm{rev}_x f, \mathrm{rev}_x g \rangle \neq \langle \mathrm{rev}_x p, \mathrm{rev}_x q \rangle  = \langle 1-x,y+1\rangle. \]

\subsection{Smith form of $S(y)$ when there are roots at infinity}

When $f$ and $g$ have infinite intersections, our proof of Theorem \ref{thm:mainresultnoinf} is not valid, and indeed we have seen in the previous subsection that the statement itself no longer holds. To circumvent this difficulty, we propose in Lemma \ref{lem:destroyinfinity} the workaround to apply a M\"{o}bius change of variable to $x$.

\begin{lemma}\label{lem:destroyinfinity}
    Let $f,g \in \FF[x,y]$ have degrees in $x$ equal to $m$ and $n$, respectively, and suppose that $\langle f,g \rangle$ is zero-dimensional. Consider the polynomials
    \[ \hat f(z,y) := (cz+d)^m f(\phi(z),y) \in \FF[z,y], \qquad  \hat g(z,y) := (cz+d)^n g(\phi(z),y) \in \FF[z,y] \]
    
    where $$x=\phi(z) := \frac{az + b}{cz + d}$$
    and $a,b,c,d \in \FF$ satisfy $ad-bc=1, c \neq 0$, and $c x_0 \neq a$
    for every $x_0 \in \FF$ such that $(x_0,y_0) \in \mathcal{V}(\langle f,g \rangle)$. Then
    \begin{enumerate}
        \item[(i)] The polynomials $\hat f$ and $\hat g$ do not have any infinite intersection of the form $(\infty,\cdot)$.
        \item[(ii)] $\langle \hat{f},\hat{g} \rangle$ is a zero-dimensional ideal of $\FF[z,y]$.
        \item[(iii)] There is a multiplicity-preserving bijection between the finite roots $(x_0,y_0)$ of $f$ and $g$ and $(\phi^{-1}(x_0),y_0)$ of $\hat f$ and $\hat g$. Furthermore, the multiplicity of the root $(\infty, y_0)$ for $f$ and $g$ is equal to the multiplicity of the root $(-d/c, y_0)$ of $\hat f$ and $\hat g$.
        \item[(iv)] The Sylvester matrices $S(y)$ of $f,g$ and $\hat S(y)$ of $\hat f,\hat g$ are strictly equivalent, and hence have the same Smith form.
        \end{enumerate}        
\end{lemma}
Observe that the existence of $a,b,c,d$ as in the statement of Lemma \ref{lem:destroyinfinity} is guaranteed because $\FF$ is an algebraically closed, and hence infinite, field\footnote{Indeed, if  $f,g \in \mathbb{K}[x,y]$ and $\mathbb{K}\subseteq \FF$ is an arbitrary infinite subfield of $\FF$, then we can always pick $\hat{f},\hat{g} \in \mathbb{K}[x,y]$ in Lemma \ref{lem:destroyinfinity}, thus ensuring that, just like $S(y)$, $\hat{S}(y)$ also has entries in $\mathbb{K}[y]$. However, this is generally not possible when $\mathbb{K}$ is a \emph{finite} subfield of $\FF$, because it may happen that there exists a point $(x_0,y_0) \in \mathcal{V}(\ideal{f,g})$ for all $x_0 \in \mathbb{K}$, thus forcing some of the coefficients $a,b,c,d$ to lie in $\FF \setminus \mathbb{K}$. Nonetheless, by Lemma \ref{lem:destroyinfinity}, we can still deduce the Smith form of $S(y)$ over $\mathbb{K}[y]$ from that of $\hat{S}(y)$ over $\FF[y]$.}. To facilitate a proof of Lemma \ref{lem:destroyinfinity}, we first prove the auxiliary Lemma \ref{lemhelp}.

\begin{lemma}\label{lemhelp}
Define $f,g,\hat{f},\hat{g},\phi,a,b,c,d$ as in Lemma \ref{lem:destroyinfinity}, and fix $x_0,y_0 \in \FF$ such that $(x_0,y_0) \in \mathcal{V}(\langle f,g\rangle)$. Let 
\[ R:=\FF[x,y]_{\langle x-x_0,y-y_0 \rangle} = \left\{  \frac{a(x,y)}{b(x,y)}, \ a(x,y),b(x,y) \in \FF[x,y], \ b(x_0,y_0)\neq 0 \right\} \] and 
\[ S:=\FF[z,y]_{\langle z-z_0,y-y_0 \rangle} = \left\{  \frac{a(z,y)}{b(z,y)}, \ a(z,y),b(z,y) \in \FF[z,y], \ b(z_0,y_0)\neq 0 \right\} \] where $z_0:=\frac{d x_0-b}{a-cx_0}$. Then,   $p(x,y) \mapsto p(\phi(z),y)$ is a ring isomorphism from $R$ to $S$. In particular,
\[ \dim_\FF R/\langle f(x,y),g(x,y) \rangle = \dim_\FF S/\langle \hat f(z,y),\hat{g}(z,y) \rangle, \]
that is, $(z_0,y_0) \in \mathcal{V}(\langle \hat{f},\hat{g}\rangle)$ and has the same multiplicity as $(x_0,y_0) \in \mathcal{V}(\ideal{f,g})$.
\end{lemma}
\begin{proof}
    We must show that the mapping is bijective. This holds because $a -c x_0 \neq 0$ by assumption, and we can then compute $c z_0 + d = (a x_0 -c)^{-1} \neq 0$. In fact, an explicit inverse is
    \[ q(z,y) \mapsto q(\phi^{-1}(x),y)= q\left(\frac{dx-b}{a-cx},y \right).     \]
This immediately implies that
\[ \dim_\FF R/\langle f(x,y),g(x,y) \rangle = \dim_\FF S/\langle  f(\phi(z),y),g(\phi(z),y) \rangle. \]
To get the statement, note in addition that $\hat f(z,y)=(cz+d)^m f(\phi(z),y)$ and $\hat g(z,y)=(cz+d)^n g(\phi(z),y)$. Hence,  since $cz+d \in S^\times$, $\hat f(z,y)$ and $f(\phi(z),y)$ are associate in $S$, and similarly for $\hat g(z,y)$ and $g(\phi(z),y)$. In particular, the two pairs generate the same ideal of $S$.   
\end{proof}

\begin{proof}[Proof of Lemma \ref{lem:destroyinfinity}]  
\begin{itemize}
    \item[(i)] Write $f=\sum_{i=0}^m f_i(y) x^i, g=\sum_{i=0}^n g_i(y) x^i$, and suppose for a contradiction that there exists $y_0 \in \FF$ such that $\hat{f}$ and $\hat{g}$ have an infinite intersection $(\infty,y_0)$. Thus, denoting by $\hat{f}_m(y),\hat{g}_n(y)$ the leading coefficients (in $z$) of $\hat{f}$ and $\hat{g}$, respectively, we have by Definition \ref{def:projectiveintersection} that $\hat{f}_m(y_0)=\hat{g}_n(y_0)=0$.

We can explicitly compute
\[  \hat{f}_m(y) = c^m \sum_{i=0}^m f_i(y) (a/c)^i = c^m f(a/c,y), \ \hat{g}_n(y) = c^n \sum_{i=0}^n g_i(y) (a/c)^i = c ^n g(a/c,y).\]

It follows $f(a/c,y_0)=g(a/c,y_0)=0$,  contradicting the assumptions.
    \item[(ii)] Suppose for a contradiction that  $h(z,y):=\gcd(\hat{f},\hat{g})$ is not a unit of $\FF[z,y]$. Denote $k=\deg_z h(z,y)$, note that $\phi^{-1}(x)=\frac{dx-b}{a-cx}$, and observe that
\[ \tilde{h}(x,y):=(a-cx)^k h\left(\phi^{-1}(x),y\right) \in \FF[x,y]\] divides both $f(x,y)=(a-cx)^m \hat{f}(\phi^{-1}(x),y)$ and $g(x,y)=(a-cx)^n \hat{g}(\phi^{-1}(x),y)$. Hence, $\tilde{h} \in \FF$, for otherwise $f$ and $g$ would have a non-unital GCD in $\FF[x,y]$; in turn, this implies $h(z,y)=\tilde{h}(cz+d)^k$. Consider now the leading coefficients of $f,g$ respectively:
\[ f_m(y) = (-c)^m \hat{f}(-d/c,y), \qquad g_n(y) = (-c)^n \hat{g}(-d/c,y).    \]
It follows that $0=h(-d/c,y)$ divides $\gcd(f_m,g_n)$ in $\FF[y]$, and hence $f_m=g_n=0$, contradicting the assumption that $m=\deg_x f$ and $n=\deg_x g$.
    \item[(iii)] Consider the notion of multiplicity given in Lemma \ref{lem:multiplicity}. The first part of the statement then follows from Lemma \ref{lemhelp}. For the second part, apply the same argument after noting  that 
     \begin{equation}\label{eq:psi}
        \mathrm{rev}_x f(x,y) = (ax-c)^m \hat{f}\left(\psi(x),y\right), \qquad \mathrm{rev}_x g(x,y) = (ax-c)^n \hat{g}\left(\psi(x),y\right)   
     \end{equation}
     where
     \[
     \psi(x):=\frac{d-bx}{ax-c} .\]
    \item[(iv)] Following \cite[Lemma 3.5]{rp3}, define the constant matrices $M_k$ as the unique elements of $\FF^{k \times k}$ such that
    \[ M_k \Lambda_k(z) = \begin{bmatrix}
        (az+b)^{k-1}\\
        (az+b)^{k-2}(cz+d)\\
        \vdots\\
        (az+b)(cz+d)^{k-2}\\
        (cz+d)^{k-1}
    \end{bmatrix} = (cz+d)^{k-1} \Lambda_k(\phi(z)) ;\]
    then, for all integers $k \geq 1$, $M_k$ is invertible by \cite[Lemma 3.5]{rp3}.

   Using Proposition \ref{prop:sproperty}, we have that

    \[  (M_n \oplus M_m) \hat{S}(y) \Lambda_{m+n}(t) = (M_n \oplus M_m) \begin{bmatrix}
        \hat{f}(t,y) \Lambda_n(t)\\
        \hat{g}(t,y) \Lambda_m(t)
    \end{bmatrix} = \begin{bmatrix}
        \hat{f}(t,y) M_n \Lambda_n(t)\\
        \hat{g}(t,y) M_m \Lambda_m(t)
    \end{bmatrix} \]
 \[ = (ct+d)^{m+n-1} \begin{bmatrix}
     f(\phi(t),y) \Lambda_n(\phi(t))\\
     g(\phi(t),y) \Lambda_m(\phi(t))
 \end{bmatrix} = S(y) (ct+d)^{m+n-1} \Lambda_{m+n}(\phi(t))  = S(y) M_{m+n} \Lambda_{m+n}(t). \]   
    Therefore, by item 1. in Proposition \ref{prop:sproperty}, $\hat{S}(y)=(M_n \oplus M_m)^{-1} S(y) M_{m+n}$, proving the statement.
\end{itemize}
\end{proof}

\begin{corollary}\label{cor414}
    Let $f,g,\hat{f},\hat{g}$ be defined as in Lemma \ref{lem:destroyinfinity}. For every pair $(x_0,z_0) \in \FF^2$ such that $x_0=\phi(z_0),z_0=\phi^{-1}(z_0)$ and for all $y_0 \in \FF$, the M\"{o}ller indices with respect to $y$ of $\langle f,g \rangle$ at $(x_0,y_0)$ coincide with the M\"{o}ller indices with respect to $y$ of $\langle \hat{f},\hat{g} \rangle$ at $(z_0,y_0)$.
\end{corollary}
\begin{proof}
    It is a direct consequence of Theorem \ref{thm:mainresultnoinf} and Lemma \ref{lem:destroyinfinity}.
\end{proof}

Theorem \ref{thm:mainresult} below generalizes Theorem \ref{thm:mainresultnoinf} to general zero-dimensional ideals.

\begin{theorem}\label{thm:mainresult}
        Let $S(y) \in \FF[y]^{(m+n) \times (m+n)}$ be the Sylvester matrix of $f,g \in \FF[x,y]$, where $m=\deg f$, $n=\deg g$. Suppose that  $ \ideal{f,g}$ is zero-dimensional.
    Fix $y_0 \in \FF$. Then:
    \begin{enumerate}
        \item If $(\infty,y_0)$ is not an infinite intersection of $f$ and $g$, then the partial multiplicities, algebraic multiplicity, and geometric multiplicity of $y_0$ as an eigenvalue of $S(y)$ are as in Theorem \ref{thm:mainresultnoinf};
        \item If $(\infty,y_0)$ is an infinite intersection of $f$ and $g$, let $x_i \in \FF$, $i=1,\dots,\ell$, denote all the elements of $\FF$ such that $f(x_i,y_0)=g(x_i,y_0)=0$, and let $\alpha_i(0),\dots,\alpha_i(\beta_i-1)$ be the M\"{o}ller indices with respect to $y$ of the the ideal $\langle f,g \rangle$ at the point $P_i=(x_i,y_0)$. Moreover, let $\alpha_\infty(0), \dots, \alpha_\infty(\beta_\infty-1)$ be the M\"{o}ller indices with respect to $y$ of $f$ and $g$ at $(\infty,y_0)$. Then, the partial multiplicities of $y_0$ in $S(y)$ are
    \[ \alpha_\infty(k), \quad k=0,\dots,\beta_\infty-1  \]
    and
\[ \alpha_{i}(k), \quad k=0,\dots,\beta_{i}-1, \quad i=1,\dots,\ell.  \]
In particular, the geometric multiplicity of $y_0$ is equal to
\[ \beta_\infty+\sum_{i=1}^\ell \beta_{i}  \]
and the algebraic multiplicity of $y_0$ is equal to
    \[  \sum_{k=0}^{\beta_\infty-1}\alpha_\infty(k)+\sum_{i=1}^\ell \sum_{k=0}^{\beta_{i}-1} \alpha_{i}(k) .\]
    \end{enumerate}
\end{theorem}
\begin{proof}
    Construct $\hat{f},\hat{g}$ and $\hat{S}(y)$ as in Lemma \ref{lem:destroyinfinity}. Since $S(y)$ and $\hat{S}(y)$ have the same invariant factors, it suffices to prove the statement for $\hat{S}(y)$, to which Theorem \ref{thm:mainresultnoinf} applies by Lemma \ref{lem:destroyinfinity}. Case 1. then follows directly by Corollary \ref{cor414}. For case 2., we must also prove that the M\"{o}ller indices with respect to $y$ of $f$ and $g$ at $(\infty,y_0)$ are equal to the M\"{o}ller indices with respect to $y$ of $\langle\hat{f},\hat{g}\rangle$ at $(a/c,y_0)$. This follows by combining Corollary \ref{cor414} and Definition \ref{def:projectiveintersection}, using also \eqref{eq:psi}.
\end{proof}

\begin{remark}\label{rem:grade}
    Observe that, even more generally, Theorem \ref{thm:mainresult} may be applied to compute the invariant factors of the $(m,n)$-Sylvester matrix where either $m>\deg_x f$ or $n > \deg_x g$ (but not both). For example, if $m=\deg_x f$ but $n>\deg_x g$, one can apply a slight variation of Theorem \ref{thm:mainresult} and Definition \ref{def:projectiveintersection}, obtained by replacing the leading coefficient in $x$ of $g$ with the zero polynomial $0 \in \FF[y]$. (Note also that the proof of item (ii) in Lemma \ref{lem:destroyinfinity} remains valid as long as either $m=\deg f$ or $n=\deg g$.)  In other words, in this scenario, in addition to all the $y$-components of the points in $\mathcal{V}(\langle f, g  \rangle)$, also every root of the leading coefficient of $f$ is an eigenvalue of $S(y)$. Their partial multiplicities correspond to the M\"{o}ller indices with respect to $y$ of $\langle \mathrm{rev}_x f, x^{n-\deg g} \mathrm{rev}_x g \rangle$. Analogous comments apply to the situation where $m>\deg f, n=\deg g$; we omit the details.
\end{remark}

\section{Smith form of B\'{e}zout matrices}\label{sec:bezoutsmith}

Theorem \ref{thm:mainresultbezout}, that we state below, is the analogue for B\'{e}zout matrices of Theorem \ref{thm:mainresult}.

\begin{theorem}\label{thm:mainresultbezout}
    Let $B(y) \in \FF[y]^{k \times k}$ be the Bezóut matrix of $f,g \in \FF[x,y]$, where $k=\max \{\deg_x f, \deg_x g\}$ and suppose $ \ideal{f,g}$ is zero-dimensional. Fix $y_0 \in \FF$. Let $P_i=(x_i,y_0)$, $i=1,\dots,\ell$, where either $x_i \in \FF$ if $(x_i,y_0) \in \mathcal{V}(I)$ or $x_i=\infty$ if $(\infty,y_0)$ is an infinite intersection of $f$ and $g$ defined with respect to the coefficients of $x^k$ in both $f$ and $g$, i.e., according to the modified version of Definition \ref{def:projectiveintersection} as described in Remark \ref{rem:grade}. Then, the partial multiplicities of $y_0$ in $B(y)$ are
    \[ \alpha_i(k), \quad k=0,\dots,\beta_i-1, \quad i=1,\dots,\ell \]
    where $\alpha_i(0),\dots,\alpha_i(\beta_i-1)$ are either the M\"{o}ller indices with respect to $y$ of $\langle f,g \rangle$ at $(x_i,y_i)$ (if $x_i \in \FF$) or (if $x_i=\infty$) the M\"{o}ller indices with respect to $y$ of $f$ and $g$ at $(\infty,y_0)$, according to the modified version of Definition \ref{def:projectiveintersection} as described in Remark \ref{rem:grade}.
\end{theorem}

Theorem \ref{thm:mainresultbezout} can be proved similarly to Theorem \ref{thm:mainresult}, and for this reason we only sketch the argument. Assume with no loss of generality (up to switching the roles of $f$ and $g$) that $k=\deg f$. The first step is to mimic Lemma \ref{lem:rootpolys} and construct root vectors for $B(y)$ corresponding to finite intersection points $(x_0,y_0)$; this is done in Lemma \ref{lem:rootvecsbez} below. If $k=\deg g$ and the leading coefficients of $f,g$ are coprime, we can then prove that the constructed set of root vectors is maximal, similarly to Lemma \ref{lem:wow} (using also that $c=1$ in Corollary \ref{cor211}). This proves Theorem \ref{thm:mainresultbezout} for the case where there are no infinite intersections; to generalize the statement and include the case of shared roots at infinity, we combine (modified versions according to Remark \ref{rem:grade} of) Definition \ref{def:projectiveintersection} and Lemma \ref{lem:destroyinfinity}.

It remains to state and prove Lemma \ref{lem:rootvecsbez}.
\begin{lemma}\label{lem:rootvecsbez}
    Let $P=(x_0,y_0) \in \mathcal{V}(\langle f, g \rangle)$, and let $\phi_i$ be the leading vector in a Gauss basis, with respect to the lex order with $x<y$, of the dual space of $\ideal{f,g}$ at $P$. Moreover, denote by $\alpha(0),\dots,\alpha(\beta-1)$ the M\"{o}ller indices with respect to $y$ of $\ideal{f,g}$ at $P$. Then, the vectors $r_i(y) = E_{(x_0,\cdot)}R_{y_0,\alpha(i)-1} (\phi_i)$ are root vectors of $B(y)$ at $y_0$, and their order is at least $\alpha(i)$.
\end{lemma}
\begin{proof}
Similarly to what happened with the Sylvester matrix, we now have $R_{y_0,\alpha} : \mathcal{D} \to \FF[x,y]^{n}$ and
\[  E_{(x_0,y_0)} D_{0,k} \Lambda_n(z)^\top B(y)  R_{y_0,\alpha}(D_{ij}) = E_{(x_0,y_0)} D_{i,j-(\alpha-k)} \BB(x,z,y) = E_{(x_0,y_0)} A^{\alpha-k}_y D_{i,j} \BB(x,z,y), \]
and for $r_i^{x_0}(y) = E_{(x_0,y)} R_{y_0,\alpha}(\phi_i)$ we have
\[ E_{(\cdot,y_0)} D_{0,k} \Lambda_n(z)^\top B(y) r_i^{x_0}(y) = E_{(x_0,y_0)} A^{\alpha-k}_y\phi_i \BB(x,z,y). \]

    Even though $\BB(x,z,y)$ is not in the ideal, for any fixed $z_0 \in \FF$, $(x-z_0)\BB(x,z_0,y)$ is. Hence by Leibnitz's rule we have
    \begin{align} \label{eq:bezleib1}
        \phi_i (x-z_0)\BB(x,z_0,y) =  (x-z_0) \phi_i \BB(x,z_0,y) +  A_x \phi_i \BB(x,z_0,y)        
    \end{align}    
    for all $\phi_i$, and in general for all $k$ we have
    \begin{align} \label{eq:bezleib2}
        A_x^k \phi_i (x-z_0) \BB(x,z_0,y) =(x-z_0) A_x^k\phi_i \BB(x,z_0,y) + A_x^{k+1} \phi_i \BB(x,z_0,y). 
    \end{align} 
    Note that the left-hand sides of \eqref{eq:bezleib1} and \eqref{eq:bezleib2} vanish at $(x_0,y_0)$, and so we may write
    \begin{align*}
         0&= E_{(x_0,y_0)}\begin{bmatrix}
\left((x-z_0)\mathcal{B}(x,z_0,y)\right)\\
\vdots\\
A^2_x \phi_i \left((x-z_0)\mathcal{B}(x,z_0,y)\right)\\
A_x \phi_i \left((x-z_0)\mathcal{B}(x,z_0,y)\right)\\
\phi_i \left((x-z_0)\mathcal{B}(x,z_0,y)\right)
\end{bmatrix} \\
&= 
\begin{bmatrix}
    x_0-z_0 & 0 &0&\dots&0 \\
    \ddots & \ddots & \ddots&&\vdots\\
    \ddots&1&x_0-z_0&0&\\
    &0&1&x_0-z_0&0\\
    0&\dots&0&1&x-z_0
\end{bmatrix}E_{(x_0,y_0)}\begin{bmatrix}
 \BB(x,z_0,y)\\
\vdots\\
 A^2_x \phi_i \BB (x,z_0,y)\\
 A_x \phi_i \BB(x,z_0,y)\\
    \phi_i \BB(x,z_0,y)
\end{bmatrix}        
    \end{align*}
    and for any $z_0 \neq x_0$ we get that $E_{(x_0,y_0)} A_x^k \phi_i \BB(x,z_0,y) = 0$ for all $k$, specifically $E_{(x_0,y_0)} \phi_i \BB(x,z_0,y) = 0$ for all $z_0 \neq x_0$. Since $\mathcal{B}(x,z,y)$ is a polynomial, this implies $E_{(x_0,y_0)} \phi_i \BB(x,z,y) = 0$ for all $z \in \FF$. 
\end{proof}

\section{Computational implications}\label{sec:numerical}
When $\FF \subseteq \mathbb{C}$, the Sylvester and B\'{e}zout polynomial matrices can be, and have been, used as a tool to numerically solve systems of $2$ polynomial equations in $2$ variables. The basic idea is as follows:  One first constructs either $S(y)$ or $B(y)$. Then, their eigenvalues can be computed via a numerical method for the polynomial eigenvalue problem \cite{MMT}, and they provide the $y$-components of the solution of the polynomial systems. Finally the $x$-component can be obtained either from the associated eigenvectors \cite{BM,GT} or by solving a univariate polynomial equation after substituting the eigenvalues $y_i$ into one of the original polynomials, say, one can solve $f(x,y_i)=0$ for each $i$ \cite{NNT}.

This approach is mathematically very elegant, but its practical implementation gives rise to various difficulties \cite{GT,NNT}. For example, if $y_0$ is a simple eigenvalue of $B(y)$ or $S(y)$, corresponding to the intersection point $(x_0,y_0) \in \mathcal{V}(\ideal{f,g})$, the condition number of $y_0$ as an eigenvalue of the resultant matrix may be as large as the square of the condition number of $(x_0,y_0)$ as a solution of $f=g=0$ \cite[Theorem 3.7 and Theorem 4.6]{NT}. While this result may appear pessimistic, this ``squared condition number" effect is a worst-case scenario that, while attained by certain specific polynomial systems \cite{NT}, does not necessarily happen for every input $f,g$ \cite{NT}; moreover, numerical algorithms exist that are able to detect and mitigate the effect of the possible increase in conditioning \cite{NNT}. Indeed, bivariate rootfinders based on either $B(y)$ or $S(y)$ have achieved remarkable success in practice  \cite{BM,GT,MD,NNT}.

\begin{remark}
    In certain applications, it is more convenient for numerical purposes to express the Sylvester and B\'{e}zout matrices in a non-monomial basis, for instance the Chebyshev basis \cite{NNT,NT}. In other words, in Definition \ref{def:bezout} one replaces  the monomial vectors $\Lambda_k$ with polynomial vectors $\Phi_k$ containing a different basis for $\FF[x]_{k-1}$; similarly, Definition \ref{def:smatrix} is changed so that Proposition \ref{prop:sproperty} is valid after replacing $\Lambda_n$ with $\Phi_n$. See \cite{NT} for more details.
    
The results obtained in this paper are still valid when the Sylvester and resultant matrices are expressed in any basis. Indeed, if $X_k \in \FF^{k \times k}$ is the change-of-basis matrix such that $\Phi_k(x)=X_k \Lambda_k(x)$, and denoting by $B(y),\hat{B}(y)$ the B\'{e}zout matrices in the monomial basis and in the basis $\Phi$, it is easy to verify that (for $k=\max \{\deg f,\deg g\}$) $\hat{B}(y)=X_k^\top B(y) X_k$, and thus $B(y)$ and $\hat{B}(y)$ have the same Smith form. Analogously, for any $m \geq \deg f$ and $n \geq \deg g$, and denoting by $\hat{S}(y),S(y)$ (resp.) the $(m,n)$-Sylvester matrix in the monomial basis and in the basis $\Phi$, by the proof of item (iv) in Lemma \ref{lem:destroyinfinity}, we have $(X_n \oplus X_m) S(y) = \hat{S}(y) X_{m+n}$, thus reaching the conclusion that the invariant factors of the Sylvester matrix do not depend on the choice of a polynomial basis.
\end{remark}

The condition number analysis of \cite{NT} applies to simple eigenvalues, i.e., algebraic multiplicity $1$. To our knowledge, no results have thus far appeared about nonsimple eigenvalues of $B(y)$ and $S(y)$. Generally, \emph{semisimple} eigenvalues (partial multiplicities all equal to $1$) of a regular polyomial matrix may be as well conditioned as simple eigenvalues, whereas \emph{defective} eigenvalues (at least one partial multiplicity higher than $1$) are ill conditioned. This fact can be made more precise \cite[Lemma 4.9]{KPM} and leads to the following rule of thumb: Employing a numerically stable algorithm \cite{MMT}, one can hope to numerically approximate the $\kappa$ eigenvalues corresponding to an exact eigenvalue $y_0$ associated with a partial multiplicity $\kappa$ up to a precision of order $\approx c {\bf u}^{1/\kappa}$, where ${\bf u}$ is the machine precision and $c$ is a constant depending on the H\"{o}lder  condition number of $y_0$ \cite{KPM} and on the employed algorithm.

As a consequence, Theorem \ref{thm:mainresult} and Theorem \ref{thm:mainresultbezout} provide clues on the achievable numerical accuracy (or lack thereof) of a resultant-based numerical approach, even in the case of non-simple eigenvalues. In particular, the M\"{o}ller indices of the ideal $\langle f,g \rangle$, or the indices at infinity of $f$ and $g$ in the case of infinite intersections, determine the partial multiplicities of each eigenvalue of the resultant matrix, and thus the first component of the H\"{o}lder condition numbers. We illustrate these observations by means of the toy Example \ref{ex:example}.
\begin{example}\label{ex:example}
Let $f=y^2+x, g=y^2-x \in \mathbb{C}[x,y]$. The ideal $\ideal{f,g}$ has dimension $2$ over $\mathbb{C}$, and the only intersection point is $O=(0,0) \in \mathbb{C}^2$, with intersection multiplicity $2$. The M\"{o}ller index with respect to $y$ of $\ideal{f,g}$ at $O$ is $2$, which implies that $0$ is a defective eigenvalue of $S(y)$ (and of $B(y)$), and thus even a stable numerical algorithm may fail to compute it accurately. On the other hand, we can switch the roles of the variables $x$ and $y$ as per Remark \ref{rem:swap}. The M\"{o}ller indices with respect to $x$ of $\ideal{f,g}$ at $O$ are $1,1$, which means that $0$ is a semisimple eigenvalue of $S(x)$ (and $B(x)$), and can be computed accurately by a stable algorithm.

For this system of equations, the M\"{o}ller indices correctly capture the fact that the $x$-component of the solution is stable with respect to small perturbations but the $y$-component is unstable. For instance, let $0 < \varepsilon \ll 5^{-1/2}$  and consider the slightly perturbed polynomial system $f+2 \varepsilon = g - \varepsilon = 0$. This corresponds to an input perturbation having norm (measured as  the Euclidean norm of the variation in  the coefficient of $f$ and $g$) equal to $\eta=\sqrt{5} \cdot \varepsilon$. A straightforward computation yields the following two complex conjugate solutions of the perturbed system:
\[ \left( - \frac{3\varepsilon}{2} ,\pm i  \sqrt{\frac{\varepsilon}{2}}\right).    \]
We conclude that the variation of the $x$-components of each solution is linear in $\eta$, whereas the variation of the $y$-components behaves as $\eta^{1/2} \gg \eta$.
\end{example}

\section*{Acknowledgements}

VN acknowledges Javier P\'{e}rez for useful discussions on some ideas that would eventually lead to this paper, as well as the University of Montana for the hospitality during a visit when those discussions happened.

\end{document}